%% file: main.tex
\documentclass[12pt,reqno]{amsart}
\usepackage{graphicx}
\usepackage{amssymb,amsmath,mathtools}
\usepackage{amsthm}
\usepackage{color,graphicx}
\usepackage{hyperref}
\usepackage{color}
\usepackage{epstopdf}
\usepackage{lpic}
\usepackage{enumerate}
\usepackage{float}
\usepackage{tikz}

\allowdisplaybreaks[4]

\newcommand{\be}{\begin{equation}}

\newcommand{\ee}{\end{equation}}
\newtheorem{theorem}{Theorem}[section]
\newtheorem{lemma}{Lemma}[section]
\newtheorem{prop}{Proposition}[section]
\newtheorem{definition}{Definition}[section]

\newtheorem{Col}{Corollary}[section]

\numberwithin{equation}{section}
\numberwithin{figure}{section}
\renewcommand{\theequation}{\thesection.\arabic{equation}}
\newtheorem{proposition}[theorem]{Proposition}

\begin{document}

\title[Simple Recursion for Mirzakhani \& Super extension]{A Simple Recursion Formula for the Mirzakhani Volume and its Super Extension}

\author{Yukun Du}
\address[Y. Du]{Department of Mathematics, University of Georgia, Athens, GA 30603}
\email{yukun.du@uga.edu}

\subjclass[2010]{Primary 14H81, 58A50}

\keywords{}

\begin{abstract}
    This paper derives a simple recursion formula for the Weil-Petersson volumes of moduli spaces of hyperbolic surfaces with boundaries. This formula demonstrates the polynomiality of the volume functions. We show that the original and our alternative formulas are equivalent by constructing their Laplace transforms. Considering these formulas' top and lowest degree terms, we recover the DVV identity and cohomology class identities for $\mathcal{M}_{g,n}$. Similar conclusions are drawn for the super-analog of these results.
\end{abstract}

\date{\today}
\maketitle
\input{Sec_1}
\input{Sec_2}
\input{Sec_3}
\input{Sec_4}
\input{Sec_5}
\input{Sec_6}
\appendix
\input{Sec_Append}
\section*{Acknowledgement}
This paper was written during my first year of Ph.D. studies at the University of California, Davis. The idea for this work was inspired by Professor Motohico Mulase, who was my initial advisor at the time. Four years later, as I was nearing the completion of my Ph.D. project, I found that my understanding of the background of this research had deepened compared to when I first wrote it. I discussed this with Prof. Mulase, and decided to revise and improve the manuscript. I am grateful to Prof. Mulase for his efforts in discussing this project with me and for his encouragement to continue exploring this area of research.

\end{document}

%% file: Sec_1.tex
\section{Main Results}\label{Chp_1}
In this paper, we explore Mirzakhani's recursion formula for the Weil-Petersson volumes of moduli spaces \cite{mirzakhani2007simple} and its super analog \cite{stanford2019jt}. Similar recursions have been previously studied, including those for Hurwitz numbers \cite{goulden1997transitive} and Catalan numbers \cite{walsh1972counting}. Such formulas are also referred to as \emph{cut-and-join equations}.

Let $\mathcal{M}_{g,n}(b_1,\dots,b_n)$ denote the moduli space of Riemann surfaces with genus $g$ and $n$ boundary components of lengths $b_1,\dots,b_n$, respectively. We consider the case where $2g+n-3\geq 0$, meaning the surfaces are hyperbolic. The hyperbolic metric on these surfaces naturally induces a K\"ahler metric on $\mathcal{M}_{g,n}(b_1,\dots,b_n)$, known as the \emph{Weil-Petersson metric}. Under this metric, the moduli space $\mathcal{M}_{g,n}(b_1,\dots,b_n)$ has a finite volume, denoted by $V_{g,n}(b_1,\dots,b_n)$. A significant result proposed by Mirzakhani states that $V_{g,n}(b_1,\dots,b_n)$ satisfies the following recursive relation:
\begin{theorem}[\cite{mirzakhani2007simple}]\label{Thm:1:1}
    The Weil-Petersson volume $V_{g,n}(\mathbf{b})$ satisfies the following recursive relation:
    \begin{equation}\label{equ:1:1}
    \begin{split}
        & \frac{\partial}{\partial b}\left(bV_{g,n}(b,\mathbf{b})\right) = \frac{1}{2}\int_0^\infty\int_0^\infty b'b''V_{g-1,n+1}(b',b'',\mathbf{b})H(b'+b'',b)db'db''\\
        & + \frac{1}{2}\sum_{\substack{g_1+g_2 = g\\ \mathbf{b}_1\sqcup \mathbf{b}_2 = \mathbf{b}}}\int_0^\infty\int_0^\infty b'b''V_{g_1,n_1}(b',\mathbf{b}_1)V_{g_2,n_2}(b'',\mathbf{b}_2)H(b'+b'',b)db'db''\\
        & + \frac{1}{2}\sum_{j=1}^{n-1}\int_0^\infty b'V_{g,n-1}(b',\mathbf{b}\backslash b_j)\left(H(b',b+b_j)+H(b',b-b_j)\right)db'.
    \end{split}
    \end{equation}
    Here, $\mathbf{b}$ denotes the tuple $(b_1,\dots,b_{n-1})$, $b$, $b'$, $b''$ and $b_j\geq 0$, and
    \[
    H(x,y) = \left(1+\exp\frac{x+y}{2}\right)^{-1}+\left(1+\exp\frac{x-y}{2}\right)^{-1}.
    \]
\end{theorem}
The construction of \eqref{equ:1:1} relies on the McShane identity. The background of the moduli space and related identities will be provided in the next section. 

It is known that Weil-Petersson volumes $V_{g,n}(\mathbf{b})$ are polynomials on the boundary lengths $\mathbf{b}$. One of the results of this paper addresses the polynomial nature of $V_{g,n}$; specifically, we derive an equivalent recursion formula for $V_{g,n}$ that does not involve the kernel $H(x,y)$:
\begin{theorem}\label{thm:1:2}
    Using the same notation as in Theorem \ref{Thm:1:1}, the volume functions $V_{g,n}(\mathbf{b})$ satisfies the following relation:
    \begin{equation}\label{equ:1:2}
    \begin{split}
        & \frac{1}{4\pi i}((b+2\pi i)V_{g,n}(b+2\pi i,\mathbf{b}) - (b-2\pi i)V_{g,n}(b-2\pi i,\mathbf{b}))\\
        & = \frac{1}{2}\iint_{\substack{b',b''\geq 0\\ b'+b''\leq b}} b'b''V_{g-1,n+1}(b',b'',\mathbf{b})db'db''\\
        & + \frac{1}{2}\sum_{\substack{g_1+g_2 = g\\ \mathbf{b}_1\sqcup \mathbf{b}_2 = \mathbf{b}}}\iint_{\substack{b',b''\geq 0\\ b'+b''\leq b}} b'b''V_{g_1,n_1}(b',\mathbf{b}_1)V_{g_2,n_2}(b'',\mathbf{b}_2)db'db''\\
        & + \frac{1}{2}\sum_{j=1}^{n-1}\left(\int_0^{b+b_j}+\int_0^{b-b_j}\right) b'V_{g,n-1}(b',\mathbf{b}\backslash b_j)db'.
    \end{split}
    \end{equation}
\end{theorem}
Utilizing \eqref{equ:1:2}, each $V_{g,n}$ can be computed simply by integrating polynomials over compact polygonal regions.

One can also consider the super analog of these hyperbolic moduli spaces. Let $\widehat{\mathcal{M}}_{g,n}(b_1,\dots,b_n)$ denote the moduli space of super Riemann surfaces with genus $n$ and $n$ boundary components, each of Neveu-Schwarz spin type and lengths $b_1,\dots, b_n$, respectively. Analogously to the definition of the Weil-Petersson volume, one defines the super volume of $\widehat{\mathcal{M}}_{g,n}(b_1,\dots,b_n)$, denoted by $\widehat{V}_{g,n}(b_1,\dots, b_n)$. By generalizing the McShane identity, Stanford and Witten extended Mirzakhani's recursion formula to the super volumes of these super moduli spaces:
\begin{theorem}[\cite{stanford2019jt}]\label{Thm:1:3}
    The super volume $\widehat{V}_{g,n}(\mathbf{b})$ satisfies the following recursive relation:
    \begin{equation}\label{equ:1:3}
    \begin{split}
        & 2\pi b\widehat{V}_{g,n}(b,\mathbf{b}) = \frac{1}{2}\int_0^\infty\int_0^\infty b'b''\widehat{V}_{g-1,n+1}(b',b'',\mathbf{b})\widehat{H}(b'+b'',b)db'db''\\
        & + \frac{1}{2}\sum_{\substack{g_1+g_2 = g\\ \mathbf{b}_1\sqcup \mathbf{b}_2 = \mathbf{b}}}\int_0^\infty\int_0^\infty b'b''\widehat{V}_{g_1,n_1}(b',\mathbf{b}_1)\widehat{V}_{g_2,n_2}(b'',\mathbf{b}_2)\widehat{H}(b'+b'',b)db'db''\\
        & + \frac{1}{2}\sum_{j=1}^{n-1}\int_0^\infty b'\widehat{V}_{g,n-1}(b',\mathbf{b}\backslash b_j)\left(\widehat{H}(b',b+b_j)+\widehat{H}(b',b-b_j)\right)db'.
    \end{split}
    \end{equation}
    Here,
    \[
    \widehat{H}(x,y) = \frac{1}{2}\left(\left(\cosh\frac{x-y}{4}\right)^{-1}-\left(\cosh\frac{x+y}{4}\right)^{-1}\right).
    \]
\end{theorem}
We adopt Norbury's notation \cite{norbury2020enumerative} for Theorem \ref{Thm:1:3}, which differs from Stanford and Witten's by constant multipliers:
\[
    \widehat{V}_{g,n}(\mathbf{b}) = (-1)^n2^{g-1+n}V^{SW}(\mathbf{b}).
\]
In this notation, the recursion formula \eqref{equ:1:3} takes a form similar to \eqref{equ:1:1} and differs from Stanford and Witten's result by certain coefficients. 

We also derive a polynomial recursion for $\widehat{V}_{g,n}(\mathbf{b})$ that is equivalent to \eqref{equ:1:3}:
\begin{theorem}\label{thm:1:4}
    Using the same notations as in Theorem \ref{Thm:1:3}, the volume functions $\widehat{V}_{g,n}(\mathbf{b})$ satisfy the following relation:
    \begin{equation}\label{equ:1:4}
    \begin{split}
        & -\frac{1}{2}\left((b+2\pi i)\widehat{V}_{g,n}(b+2\pi i,\mathbf{b}) + (b-2\pi i)\widehat{V}_{g,n}(b-2\pi i,\mathbf{b})\right)\\
        & = \frac{1}{2}\int_0^b b'(b-b')\widehat{V}_{g-1,n+1}(b',b-b',\mathbf{b})db'\\
        & + \frac{1}{2}\sum_{\substack{g_1+g_2 = g\\ \mathbf{b}_1\sqcup \mathbf{b}_2 = \mathbf{b}}}\int_0^b b'(b-b')\widehat{V}_{g_1,n_1}(b',\mathbf{b}_1)\widehat{V}_{g_2,n_2}(b-b',\mathbf{b}_2)db'\\
        & + \frac{1}{2}\sum_{j=1}^{n-1}\left((b+b_j)\widehat{V}_{g,n-1}(b+b_j,\mathbf{b}\backslash b_j)+(b-b_j)\widehat{V}_{g,n-1}(b-b_j,\mathbf{b}\backslash b_j)\right).
    \end{split}
    \end{equation}
\end{theorem}
Do and Norbury suggested other relations for $V_{g,n}$ and $\widehat{V}_{g,n}$ involving a complex boundary length $2\pi \mathrm{i}$ and a multiplier $(2g-3+n)$:
\begin{prop}[\cite{do2009weil}]\label{prop:1:1}
    For any $g,n$ such that $2g+n\geq 3$,
    \[
        \frac{\partial}{\partial b}V_{g,n}(2\pi\mathrm{i},\mathbf{b}) = 2\pi\mathrm{i}(2g-3+n)V_{g,n-1}(\mathbf{b}).
    \]
    Here, $\frac{\partial}{\partial b}$ differentiates on the first variable of $V_{g,n}$.
\end{prop}
\begin{prop}[\cite{norbury2020enumerative}]\label{prop:1:2}
    For any $g,n\geq 1$,
    \[
        \widehat{V}_{g,n}(2\pi\mathrm{i},\mathbf{b}) = (2g-3+n)V_{g,n-1}(\mathbf{b}).
    \]
\end{prop}

\emph{Topological recursions} are established by applying Laplace transformations to cut-and-join equations, as discussed in \cite{dumitrescu2015lectures,mulase2012spectral}. In this paper, we derive topological recursions from the two cut-and-join equations for the Weil-Petersson volume $V_{g,n}(\mathbf{b})$.
\begin{theorem}\label{thm:1:5}
    Let $F^V_{g,n}(\mathbf{t})$ denote the Laplace transform of the Weil-Petersson volume $V_{g,n}(\mathbf{b})$, i.e.,
    \[
    F_{g,n}^V(\mathbf{t}) = \int_{\mathbb{R}_+^n}V_{g,n}(\mathbf{b})e^{-\mathbf{b}\cdot\mathbf{t}}d\mathbf{b}.
    \]
    Then,
    \begin{enumerate}
        \item Equation \eqref{equ:1:1} implies the following topological recursion for $F_{g,n}^V$:
        \begin{equation}\label{equ:1:5}
        \begin{split}
            & -t\frac{\partial}{\partial t}F_{g,n}^V(t,\mathbf{t}) = \mathsf{P}_t\left(\pi\csc(2\pi t)\left.\frac{\partial^2}{\partial t\partial t'}F_{g-1,n+1}^V(t,t',\mathbf{t})\right|_{t' = t}\right) \\
            & + \sum_{\substack{g_1+g_2 = g\\\mathbf{t}_1\sqcup \mathbf{t}_2 = \mathbf{t}}}\mathsf{P}_t\left(\pi\csc(2\pi t)\left(\frac{\partial}{\partial t}F^V_{g_1,n_1}(t,\mathbf{t}_1)\right)\left(\frac{\partial}{\partial t}F^V_{g_2,n_2}(t,\mathbf{t}_2)\right)\right) \\
            & + \sum_{j=1}^{n-1}\mathsf{P}_t\left(2\pi t_j(t^2 - t_j^2)^{-1}\csc(2\pi t)\frac{\partial}{\partial t}F_{g,n-1}^V(t,\mathbf{t}\backslash t_j)\right).
        \end{split}
        \end{equation}
        Here, $\mathbf{t}$ is a tuple $(t_1,\dots,t_{n-1})$, and $\mathsf{P}_t$ denotes taking the principal part of a Laurent series at $t = 0$.
        \item Equation \eqref{equ:1:2} implies the following topological recursion for $F_{g,n}^V$:
        \begin{equation}\label{equ:1:6}
        \begin{split}
            & -\mathsf{P}_t\left(\frac{\sin2\pi t}{2\pi}\frac{\partial}{\partial t}F_{g,n}^V(t,\mathbf{t})\right) = \frac{1}{2t}\left.\frac{\partial^2}{\partial t\partial t'}F_{g-1,n+1}^V(t,t',\mathbf{t})\right|_{t' = t}\\
            & + \sum_{\substack{g_1+g_2 = g\\\mathbf{t}_1\sqcup \mathbf{t}_2 = \mathbf{t}}} \frac{1}{2t}\left(\frac{\partial}{\partial t}F^V_{g_1,n_1}(t,\mathbf{t}_1)\right)\left(\frac{\partial}{\partial t}F^V_{g_2,n_2}(t,\mathbf{t}_2)\right)\\
            & + \sum_{j=1}^{n-1}(t^2 - t_j^2)^{-1}\left(\frac{t_j}{t}\frac{\partial}{\partial t}F_{g,n-1}^V(t,\mathbf{t}\backslash t_j) - \frac{t}{t_j}\frac{\partial}{\partial t_j}F_{g,n-1}^V(\mathbf{t})\right).
        \end{split}
        \end{equation}
    \end{enumerate}
    Moreover, equations \eqref{equ:1:5} and \eqref{equ:1:6} are equivalent.
\end{theorem}
Eynard and Orantin considered a different form of Laplace transform in \cite{eynard2007weil}:
\[
F_{g,n}^{EO}(\mathbf{t}) = \frac{(-1)^n\partial^n}{\partial t_1\dots\partial t_n}F_{g,n}^V(\mathbf{t}).
\]
Using this notation, they derived a topological recursion formula for the Weil-Petersson volume. In Section \ref{Chp_3}, we proved that Eynard and Orantin's formula is equivalent to equation \eqref{equ:1:5}. 

We also obtain the super analog of Theorem \ref{thm:1:5}:
\begin{theorem}\label{thm:1:6}
    Let $F^{\widehat{V}}_{g,n}(\mathbf{t})$ denote the Laplace transform of the super volume $\widehat{V}_{g,n}(\mathbf{b})$, i.e.,
    \[
    F_{g,n}^{\widehat{V}}(\mathbf{t}) = \int_{\mathbb{R}_+^n}\widehat{V}_{g,n}(\mathbf{b})e^{-\mathbf{b}\cdot\mathbf{t}}d\mathbf{b}.
    \]
    Then,
    \begin{enumerate}
        \item Equation \eqref{equ:1:3} implies the following topological recursion for $F_{g,n}^{\widehat{V}}$:
        \begin{equation}\label{equ:1:7}
            \begin{split}
                & -\frac{\partial}{\partial t}F_{g,n}^{\widehat{V}}(t,\mathbf{t}) = \mathsf{P}_t\left(\frac{1}{2}\sec(2\pi t)\left.\frac{\partial^2}{\partial t\partial t'}F_{g-1,n+1}^{\widehat{V}}(t,t',\mathbf{t})\right|_{t' = t}\right) \\
                & + \sum_{\substack{g_1+g_2 = g\\\mathbf{t}_1\sqcup \mathbf{t}_2 = \mathbf{t}}}\mathsf{P}_t\left(\frac{1}{2}\sec(2\pi t)\left(\frac{\partial}{\partial t}F^{\widehat{V}}_{g_1,n_1}(t,\mathbf{t}_1)\right)\left(\frac{\partial}{\partial t}F^{\widehat{V}}_{g_2,n_2}(t,\mathbf{t}_2)\right)\right) \\
                & + \sum_{j=1}^{n-1}\mathsf{P}_t\left(t_j(t^2 - t_j^2)^{-1}\sec(2\pi t)\frac{\partial}{\partial t}F_{g,n-1}^{\widehat{V}}(t,\mathbf{t}\backslash t_j)\right).
            \end{split}
        \end{equation}
        \item Equation \eqref{equ:1:4} implies the following topological recursion for $F_{g,n}^{\widehat{V}}$:
        \begin{equation}\label{equ:1:8}
            \begin{split}
                & -\mathsf{P}_t\left(\cos(2\pi t)\frac{\partial}{\partial t}F_{g,n}^{\widehat{V}}(t,\mathbf{t})\right) = \frac{1}{2}\left.\frac{\partial^2}{\partial t\partial t'}F_{g-1,n+1}^{\widehat{V}}(t,t',\mathbf{t})\right|_{t' = t}\\
                & + \frac{1}{2}\sum_{\substack{g_1+g_2 = g\\\mathbf{t}_1\sqcup \mathbf{t}_2 = \mathbf{t}}} \left(\frac{\partial}{\partial t}F^{\widehat{V}}_{g_1,n_1}(t,\mathbf{t}_1)\right)\left(\frac{\partial}{\partial t}F^{\widehat{V}}_{g_2,n_2}(t,\mathbf{t}_2)\right)\\
                & + \sum_{j=1}^{n-1}t_j(t^2 - t_j^2)^{-1}\left(\frac{\partial}{\partial t}F_{g,n-1}^{\widehat{V}}(t,\mathbf{t}\backslash t_j) - \frac{\partial}{\partial t_j}F_{g,n-1}^{\widehat{V}}(\mathbf{t})\right).
            \end{split}
        \end{equation}
    \end{enumerate}
    Moreover, equations \eqref{equ:1:7} and \eqref{equ:1:8} are equivalent.
\end{theorem}
Equation \ref{equ:1:7} is also appeared in \cite{norbury2020enumerative}.

Topological recursion formulas contain a wealth of information. For instance, by examining the highest- and lowest-degree terms of the Hurwitz topological recursion, one recovers the Witten-Konsevich theorem and the $\lambda_g$ formula, respectively \cite{mulase2009polynomial}. In this paper, we analyze the top- and lowest-degree terms of the new cut-and-join equations \eqref{equ:1:2} and \eqref{equ:1:4}, recovering the Virasoro constraint equations as well as equations of cohomology classes, respectively.

The paper is organized as follows. In Section \ref{Chp_1.5}, we review the background material, including the moduli space of hyperbolic Riemannian surfaces and its Weil-Petersson volume, Mirzakhani's recursion formula, and their super analogs. In Section \ref{Chp_2}, we provide proofs of our newly proposed recursion relations \eqref{equ:1:2} and \eqref{equ:1:4} and offer a new perspective on the polynomiality of the Weil-Petersson volume functions. In Section \ref{Chp_3}, by considering the Laplace transformation, we derive topological recursions for both the Weil-Petersson volume and the super volume, revealing additional relationships between our recursion formulas \eqref{equ:1:2} and \eqref{equ:1:4} and the original ones. In Section \ref{Chp_4}, we analyze the highest- and lowest-degree terms of the topological recursions to recover specific equations analogously to the Hurwitz number case.

%% file: Sec_2.tex
\section{Background}\label{Chp_1.5}
\subsection{The Volume of Hyperbolic Moduli Spaces}
For any $g,n\in \mathbb{N}$ with $2g+n\geq 3$, define the moduli space $\mathcal{M}_{g,n}(b_1,\dots,b_n)$ as the space consisting of all hyperbolic structures on the topological surface of genus $g$ with $n$ geodesic boundary components, each of length $b_i$. More rigorously,
\[
\mathcal{M}_{g,n}(\mathbf{b}) = \left\{M\left|\begin{array}{l}
     M\ \mathrm{is\ a\ hyperbolic\ surface\ of\ genus\ } g\mathrm{\ with} \\ \mathrm{geodesic\ boundary\ components\ of\ length}\ b_i \end{array}\right\}\right/\sim,
\]
where $\mathbf{b} = (b_1,\dots,b_n)$ and $M\sim M'$ if and only if there exists an isometry taking $M$ and its boundary $\partial M$ to $M'$ and $\partial M'$, respectively.

To better characterize such hyperbolic surfaces in $\mathcal{M}_{g,n}(\mathbf{b})$, one considers \emph{pants decompositions}:
\begin{definition}
    A \textbf{pants decomposition} of a hyperbolic surface $M$ of genus $g$ and $n$ boundary components is a collection $\{\alpha_i\}_{i=1}^{3g-3+n}$ of simply closed geodesic curves in $\mathrm{int}(M)$, such that $M\backslash\left(\bigcup_{i=1}^{3g-3+n}\alpha_i\right)$ consists of $(2g-2+n)$ connected components, each of which is a thrice-punctured sphere, also known as a ``pair of pants''.
\end{definition}
We note that the number of selected simply closed geodesic curves and which of pairs of pants are uniquely determined by the topological nature of the surface $M$. For a pair of pants with specified boundary component lengths, there is a unique hyperbolic structure on it:
\begin{proposition}[\cite{thurston1997three}]
    For any $b_1,b_2,b_3>0$, a pair of pants with boundary components of lengths $2b_1$, $2b_2$, and $2b_3$ can be obtained by gluing together two copies of a right-angled hyperbolic hexagon with side lengths $b_1$, $\beta_1$, $b_2$, $\beta_2$, $b_3$, and $\beta_3$. Here, $\beta_i$, for $i=1,2,3$, are given by
    \[
    \cosh \beta_i = \frac{\cosh b_i + \cosh b_{i+1}\cosh b_{i+2}}{\sinh b_{i+1}\sinh b_{i+2}},
    \]
    with subscripts taken modulo $3$.
\end{proposition}
On the other hand, there are various ways to assemble a given collection of $(2g-2+n)$ pairs of pants along $(3g-3+n)$ pairs of identified boundary components. Specifically, suppose that $\alpha_i$ and $\alpha_i'$ are identified boundary components. If we regard them as unit speed curves $\alpha_i(t)$ and $\alpha_i'(t)$, $t\in\mathbb{R}$. Then a possible identification would be
\[
\alpha_i(t)\sim \alpha_i'(t+\tau_i),\ \forall t\in \mathbb{R},
\]
where $\tau_i\in\mathbb{R}$ is called a \emph{twist parameter} for $\alpha_i$ and $\alpha_i'$. From the above discussion, it is evident that lengths $l_i\in\mathbb{R}_+$ and twist parameters $\tau_i\in\mathbb{R}$ for closed geodesics $\alpha_i$, $i=1,\dots, 3g-3+n$, uniquely determine a hyperbolic structure on $M$.

However, a hyperbolic surface admits different pants decompositions with varying geodesic lengths and twist parameters. Therefore, we consider to associate a hyperbolic structure in $\mathcal{M}_{g,n}(\mathbf{b})$ with a specified pants decomposition; that is, we fix a smooth surface $\Sigma_{g,n}$ of genus $g$ and $n$ boundary components and assign $M$ with a diffeomorphism from $\Sigma_{g,n}$, resulting in a larger space.
\begin{definition}
    The \textbf{Teichm\"uller space} for hyperbolic surfaces of genus $g$ and $n$ boundary components of lengths $b_1,\dots,b_n$ is
    \[
    \mathcal{T}_{g,n}(\mathbf{b}) = \left\{(M,f)\left|\begin{array}{l}
     M\ \mathrm{is\ a\ hyperbolic\ surface\ of\ genus\ } g\mathrm{\ with \ geodesic\ boundary}\\ \mathrm{components\ of\ length}\ b_i,\ f:\Sigma_{g,n}\to M\ \mathrm{is\ a\ diffeomorphism} \end{array}\right\}\right/\sim,
    \]
    where $(M,f)\sim (M',f')$ if and only if there exists an isometry $\varphi: M\to M'$ such that $\varphi\circ f$ is isotopic to $f'$.
\end{definition}
By assigning $\Sigma_{g,n}$ a (topological) pants decomposition, the diffeomorphism $f: \Sigma_{g,n}\to M$ specify a pants decomposition for $M$. In fact:
\begin{theorem}
    The map $\mathcal{T}_{g,n}(\mathbf{b})\to \mathbb{R}_+^{3g-3+n}\times\mathbb{R}^{3g-3+n}$ that takes a pair $(M,f)$ to the corresponding geodesic lengths and twist parameters is a homeomorphism.
\end{theorem}
We refer to the collection $(l_1,\dots, l_{3g-3+n},\tau_1,\dots, \tau_{3g-3+n})\in \mathbb{R}_+^{3g-3+n}\times\mathbb{R}^{3g-3+n}$ as the \emph{Fenchel-Nielsen coordinate} for $\mathcal{T}_{g,n}(\mathbf{b})$, \cite{wolpert1982fenchel}.

The moduli space $\mathcal{M}_{g,n}(\mathbf{b})$ is realized as a quotient of the Teichm\"uller space:
\[
\mathcal{M}_{g,n}(\mathbf{b}) = \mathcal{T}_{g,n}(\mathbf{b})/\mathrm{Mod}_{g,n},
\]
where $\mathrm{Mod}_{g,n} = \mathrm{Diff}^+(\Sigma_{g,n})/\mathrm{Diff}_0(\Sigma_{g,n})$, known as the \emph{mapping class group}, is a discrete group consisting of the orientation-preserving diffeomorphism group modulo the identity component of the diffeomorphism group of $\Sigma_{g,n}$.

The hyperbolic structures on the surfaces in the Teichm\"uller space $\mathcal{T}_{g,n}(\mathbf{b})$ induce a metric on the Teichm\"uller space as follows. The Kodaira-Spencer identification\cite{kodaira1986complex} implies that
\[
T_M\mathcal{T}_{g,n}(\mathbf{b}) = H^1(M,\mathcal{T}_M),
\]
where $\mathcal{T}_M$ is the sheaf of holomorphic vector fields on $M$. Since $H^2(M,\mathcal{T}_M) = \mathcal{T}_M$, one defines a bilinear form:
\[
\langle u,v\rangle = \int_M uv,\ \forall u,v\in T_M\mathcal{T}_{g,n}(\mathbf{b}) = H^1(M,\mathcal{T}_M).
\]
\begin{proposition}[\cite{wolpert1981elementary}]
    The aforementioned inner product on $T_M\mathcal{T}_{g,n}(\mathbf{b})$ defines a K\"ahler metric, called the \textbf{Weil-Petersson metric}. It is invariant under the action of $\mathrm{Mod}_{g,n}$, and is expressed in the Fenchel-Nielsen coordinate as
    \[
    \omega = \sum_{i=1}^{3g-3+n}dl_i\wedge d\tau_i.
    \]
\end{proposition}
The Weil-Petersson metric induces a volume form on the Teichm\"uller space, which is also invariant under the $\mathrm{Mod}_{g,n}$-action. Cosequently, one considers the volume of the moduli space $\mathcal{M}_{g,n}(\mathbf{b}) = \mathcal{T}_{g,n}(\mathbf{b})/\mathrm{Mod}_{g,n}$, denoted by $V_{g,n}(\mathbf{b})$. For convenience, we define $V_{1,1}(b)$ as half of the volume of $\mathcal{M}_{1,1}(b)$, since the surfaces in $\mathcal{M}_{1,1}(b)$ have non-trivial self-isometries of order two, \cite{witten1990two}. Computations yield
\begin{equation}\label{equ:2:1}
    V_{0,3}(\mathbf{b}) = 1,\ V_{1,1}(b) = \frac{1}{48}(b^2+4\pi^2).
\end{equation}

The Weil-Petersson volumes $V_{g,n}(\mathbf{b})$ can be determined using a recursion formula, derived from a generalized McShane identity, \cite{mirzakhani2007simple}:
\begin{theorem}
    Let $M$ be a hyperbolic surface with genus $g$, let $\beta$ be a geodesic boundary component of length $b$, with other boundary components be $\beta_i$ of length $b_i$, $i=1,\dots,n-1$. Then,
    \begin{equation}\label{Equ_McSh}
    b = \sum_{(\alpha_j,\alpha_k)}\mathcal{D}(b,l(\alpha_j),l(\alpha_k))+\sum_{i=1}^{n-1}\sum_{\alpha_j}\mathcal{R}(b,b_i,l(\alpha_j)),
    \end{equation}
    where functions $\mathcal{D}(x,y,z)$ and $\mathcal{R}(x,y,z)$ are defined as follows:
    \[
    \begin{split}
        &\mathcal{D}(0,y,z) = 0,\ \frac{\partial}{\partial x}\mathcal{D}(x,y,z) = H(y+z,x)\\
        &\mathcal{R}(0,y,z) = 0,\ \frac{\partial}{\partial x}\mathcal{R}(x,y,z) = \frac{1}{2}(H(z,x+y)+H(z,x-y)),
    \end{split}
    \]
    and $H(x,y) = \left(1+\exp\frac{x+y}{2}\right)^{-1}+\left(1+\exp\frac{x-y}{2}\right)^{-1}$. The first summation is over unordered pairs of simply closed geodesics in $M$ that bound a pair of pants with $\beta$, and the second summation is over simple closed geodesics in $M$ that bound a pair of pants with $\beta$ and $\beta_i$.
\end{theorem}
\begin{proof}[Outline of the proof]
    For each point $X\in \beta$, denote by $\gamma_X:\mathbb{R}_{\geq 0}\to M$ the geodesic on $M$ that intersects $\beta$ orthogonally at $X = \gamma_X(0)$. As $t$ increases, one of the following cases holds for $\gamma_X$:
    \begin{enumerate}
        \item The geodesic $\gamma_X$ neither intersects itself nor reaches a boundary component.
        \item The geodesic $\gamma_X$ intersects itself or the boundary component $\beta$.
        \item The geodesic $\gamma_X$ reaches the boundary component $\beta_i$ for $i=1,\dots,n-1$.
    \end{enumerate}
    Mirzakhani demonstrated the following, \cite{mirzakhani2007simple}:
    \begin{enumerate}
        \item The points $X$ in case (1) form a Cantor set.
        \item The points $X$ in case (2) form a collection of open intervals $(X_{jk},X_{jk}')$, such that $\gamma_{X_{jk}}$ and $\gamma_{X_{jk}'}$ are asymptotic to simply closed geodesics $\alpha_j$ and $\alpha_k$, respectively. For any $X\in (X_{jk},X_{jk}')$, a sufficiently small neighborhood of $\beta\cup \gamma_X$ is topologically a pair of pants, with its other two boundary components homotopic to $\alpha_j$ and $\alpha_k$.
        \item The points $X$ in case (3) form a collection of open intervals $(X_{j},X_{j}')$, such that both $\gamma_{X_{j}}$ and $\gamma_{X_{j}'}$ are asymptotic to the simply closed geodesic $\alpha_j$. Furthermore, for any $X\in (X_{j},X_{j}')$, a sufficiently small neighborhood of $\beta\cup\beta_i\cup \gamma_X$ is topologically a pair of pants, with its third boundary component homotopic to $\alpha_j$.
    \end{enumerate}
    The length of a particular interval $(X_{jk},X_{jk}')$ or $(X_{j},X_{j}')$ can be computed by considering only the corresponding pair of pants (rather than the entire surface $M$) using elementary hyperbolic trigonometry:
    \begin{enumerate}
        \item The length of an interval $(X_{jk},X_{jk}')$ in case (2) is one half of $\mathcal{D}(b,l(\alpha_j),l(\alpha_k))$.
        \item The length of an interval $(X_{j},X_{j}')$ in case (3) with respect to $\beta_i$ is $\mathcal{R}(b,b_i,l(\alpha_j))$.
    \end{enumerate}
    Moreover, there are exactly two such intervals for each pair $(\alpha_j,\alpha_k)$ in case (2), and exactly one such interval for each $\alpha_j$ in case (3) with respect to each $\beta_i$. Summing up the lengths of all these intervals yields the identity.
\end{proof}
This identity generalizes the original McShane identity \cite{mcshane1998simple} for the special case of $\mathcal{M}_{1,1}(0)$. To derive Mirzakhani's recursive formula \eqref{equ:1:1}, we further categorize the intervals in case (2) into two subtypes:
\begin{itemize}
    \item The complement of the pair of pants bounded by $\beta$, $\alpha_j$ and $\alpha_k$ is connected (thus forming a surface of genus $(g-1)$ and $(n+1)$ boundary components).
    \item The complement of the pair of pants bounded by $\beta$, $\alpha_j$ and $\alpha_k$ consists of two connected components of genera $g_1$ and $g_2$, with $g_1 + g_2 = g$. Moreover, the lengths of their other boundary components constitute tuples $\mathbf{b}_1$ and $\mathbf{b}_2$, such that $\mathbf{b}_1\sqcup \mathbf{b}_2 = (b_1,\dots,b_{n-1})$.
\end{itemize}
For any $M\in\mathcal{M}_{g,n}(b,\mathbf{b})$, denote the sum of interval lengths in these subtypes by $\mathcal{D}_{con}(M)$ and $\mathcal{D}_{g_1,\mathbf{b_1}}(M)$, respectively. Additionally, denote the sum of the interval lengths of type (3) with respect to $\beta_i$ by $\mathcal{R}_i(M)$, $i=1,\dots, n-1$. Integrating \eqref{Equ_McSh} over $\mathcal{M}_{g,n}(\mathbf{b})$ gives:
\[
bV_{g,n}(b,\mathbf{b}) = \int_{\mathcal{M}_{g,n}}\mathcal{D}_{con}(M) + \sum_{\substack{g_1+g_2 = g\\\mathbf{b}_1\sqcup\mathbf{b}_2 = \mathbf{b}}}\int_{\mathcal{M}_{g,n}}\mathcal{D}_{g_1,\mathbf{b_1}}(M) + \sum_{i=1}^{n-1}\int_{\mathcal{M}_{g,n}}\mathcal{R}_i(M).
\]
This ultimately leads to \eqref{equ:1:1} by changing the domain of integrals from $\mathcal{M}_{g,n}(b,\mathbf{b})$ to its universal covering $\mathcal{T}_{g,n}(b,\mathbf{b})$.

The polynomiality of the Weil-Petersson volumes $V_{g,n}(\mathbf{b})$ is natural. As shown in \cite{mirzakhani2007weil}, the Weil-Petersson form $\omega$ naturally extends to the boundary, yielding a compactification $\overline{\mathcal{M}}_{g,n}$ of $V_{g,n}(\mathbf{b})$. The volume of $V_{g,n}(\mathbf{b})$ is then expressed as
\begin{equation}\label{is:1}
V_{g,n}(\mathbf{b}) = \int_{\overline{\mathcal{M}}_{g,n}}\exp\left(2\pi^2\kappa_1 + \frac{1}{2}\sum b_i^2\psi_i\right),
\end{equation}
where $\psi_i$ is the Chern class over the $i$-th tautological line bundle, and $\kappa_1 = [\omega]/2\pi$ is the \emph{first Mumford tautological class}, \cite{wolpert1983homology}. It follows that $V_{g,n}(\mathbf{b})$ is a even polynomial over $b_1,\dots,b_n$, where the coefficients are integrals of products of Chern classes and the Mumford class.
\subsection{Super Hyperbolic Geometry and Super Moduli Spaces}
Stanford and Witten, \cite{stanford2019jt} generalized Mirzakhani's work to the moduli space of super hyperbolic surfaces, \cite{hodgkin1988super,crane1988super}. We begin by introducing super numbers:
\begin{definition}
    The algebra of (real) \textbf{super numbers} is a $(\mathbb{Z}/2\mathbb{Z})$-graded algebra:
    \[
    \widehat{\mathbb{R}} = \left\{\left.a = a_{\#} + \sum_{k=1}^n\sum_{i_1<\dots<i_k}a_{i_1\dots i_k}e_{i_1}\dots e_{i_k}\right|n,i_1,\dots,i_k\in \mathbb{N}_+,\ a_{\#}, a_{i_1\dots i_k}\in\mathbb{R}\right\},
    \]
    where $e_1$, $e_2$,\dots are generators that anti-commute with each other. The term $a_\#$ of degree zero is called the \textbf{body} of $a$. The grading is given as $\widehat{\mathbb{R}} = \widehat{\mathbb{R}}[0]\otimes \widehat{\mathbb{R}}[1]$, consisting of super numbers with even- and odd-degree terms, respectively. They are also known as the \textbf{Bosonic} and \textbf{Fermionic} parts of super numbers.

    A super number $a$ is \textbf{invertible} if its body $a_\# \neq 0$. We say that two super numbers satisfy $a\geq b$ if their bodies satisfy $a_\#\geq b_\#$ as real numbers.

    The super vector space over $\widehat{\mathbb{R}}$ of dimension $m|n$ is defined as
    \[
    \widehat{\mathbb{R}}^{m|n} = \{(x_1,\dots,x_m|\varphi_1,\dots,\varphi_n)|x_i\in\widehat{\mathbb{R}}[0],\ \varphi_i\in \widehat{\mathbb{R}}[1]\}.
    \]
    The group of endomorphisms of $\widehat{\mathbb{R}}^{m|n}$ is
    \[
    M_{m|n}(\widehat{\mathbb{R}}) = \left\{\left.\begin{pmatrix}
        A & B\\ C & D
    \end{pmatrix}\right|\begin{array}{c}
        A\in Mat_{m\times m}(\widehat{\mathbb{R}}[0]),\ D\in Mat_{n\times n}(\widehat{\mathbb{R}}[0]),\\
        B\in Mat_{m\times n}(\widehat{\mathbb{R}}[1]),\ C\in Mat_{n\times m}(\widehat{\mathbb{R}}[1]) 
    \end{array}\right\}.
    \]
    The \textbf{Berezinian} of an element in $M_{m|n}(\widehat{\mathbb{R}})$ is defined as
    \[
    \mathrm{Ber}\begin{pmatrix}
        A & B\\ C & D
    \end{pmatrix} = \frac{\det(A-BD^{-1}C)}{\det(D)},
    \]
    generalizing the concept of determinants. The \textbf{supertranspose} of an element in $M_{m|n}(\widehat{\mathbb{R}})$ is defined as
    \[
    \begin{pmatrix}
        A & B\\ C & D
    \end{pmatrix}^{\mathrm{ST}} = \begin{pmatrix}
        A^\mathrm{T} & C^\mathrm{T}\\ -B^\mathrm{T} & D^\mathrm{T}
    \end{pmatrix}.
    \]
\end{definition}
We use a generalized hyperboloid model to define the super hyperbolic plane:
\begin{definition}
    Equip the super vector space $\widehat{\mathbb{R}}^{2,1|2}$ with the bilinear form
    \[
    \langle(x,y,z|\varphi,\psi),(x',y',z'|\varphi',\psi')\rangle = xx' - yy' - zz' + \varphi\psi' + \varphi'\psi.
    \]
    The \textbf{super hyperbolic plane} is defined as a hypersurface in $\widehat{\mathbb{R}}^{2,1|2}$:
    \[
    \widehat{\mathbf{H}}^2 = \{X\in \widehat{\mathbb{R}}^{2,1|2}|\langle X,X\rangle = 1\},
    \]
    endowed with the (Bosonic-valued) Riemannian supermetric
    \[
    d(X_1,X_2) = \mathrm{arccosh}\langle X_1,X_2\rangle.
    \]
\end{definition}
Similarly to the $SL(2,\mathbb{R})$-action on $\mathbf{H}^2$, one considers the following group action of isometries on $\widehat{\mathbf{H}}^2$:
\begin{definition}
    Define the group $OSp(1|2,\widehat{\mathbb{R}})$ as
    \[
    OSp(1|2,\widehat{\mathbb{R}}) = \{g\in M_{2|1}(\widehat{\mathbb{R}})|g^{\mathrm{ST}}Jg = J,\ \mathrm{Ber}(g) = 1\},
    \]
    where
    \[
    J = \begin{pmatrix}
        0 & 1 & 0\\ -1 & 0 & 0\\ 0 & 0 & -1
    \end{pmatrix}.
    \]
    The group action $OSp(1|2,\widehat{\mathbb{R}})\curvearrowright\widehat{\mathbf{H}}^2$ is given by
    \[
    g.X = g^{\mathrm{ST}}Xg,
    \]
    where $X = (x,y,z|\varphi,\psi)\in \widehat{\mathbb{R}}^{2,1|2}$ is realized as the matrix
    \[
    \begin{pmatrix}
        x+z & y & \varphi \\ y & x-z & \psi \\ -\varphi & -\psi & 0
    \end{pmatrix}.
    \]
\end{definition}
Next, we consider the moduli space of super hyperbolic surfaces. The group $OSp(1|2,\widehat{\mathbb{R}})$ has a canonical decomposition, \cite{penner2021super}:
\[
\begin{pmatrix}
    a & b & 0\\ c & d & 0\\ 0 & 0 & 1
\end{pmatrix}\cdot \begin{pmatrix}
    1-\alpha\beta/2 & 0 & \alpha\\ 0 & 1-\alpha\beta/2 & \beta \\ \beta & -\alpha & 1+\alpha\beta
\end{pmatrix}.
\]
The image of a geodesic in $\widehat{\mathbf{H}}$ under the Fermionic part of the decomposition provides a \emph{thickening} of the geodesic, which is of dimension $(1|2)$; these gives locally the geodesic boundaries of super hyperbolic surfaces\cite{stanford2019jt}. 

However, a super pair of pants with given boundary component lengths has additional degrees of freedom. These are captured by considering the holonomies around the three holes of a pair of pants: \cite{stanford2019jt}
\[
\gamma_a = \delta_a\begin{pmatrix}
    e^{a/2} & * & 0 \\ 0 & e^{-a/2} & 0\\ 0 & 0 & \delta_a
\end{pmatrix}\exp\begin{pmatrix}
    0 & 0 & \xi \\ 0 & 0 & 0\\ 0 & -\xi & 0
\end{pmatrix},\ \gamma_b = \delta_b\begin{pmatrix}
    e^{b/2} & 0 & 0 \\ 1 & e^{-b/2} & 0\\ 0 & 0 & \delta_b
\end{pmatrix}\exp\begin{pmatrix}
    0 & 0 & 0 \\ 0 & 0 & \psi\\ \psi & 0 & 0
\end{pmatrix},
\]
with the condition $\gamma_a\gamma_b\gamma_c = Id$, and $\gamma_c$ is conjugate to $\delta_c diag(e^{c/2},e^{-c/2},\delta_c)$. Here, $\psi$ and $\xi$ are fermionic moduli, and $\delta_a$, $\delta_b$ and $\delta_c\in\{-1,1\}$, satisfying $\delta_a\delta_b\delta_c = -1$. The boundary $a$ (similarly $b$ or $c$) is of the Ramond (R) spin type if $\delta_a = 1$, and of Neveu-Schwarz (NS) spin type if $\delta_a = -1$. We thus define:
\[
\widehat{\mathcal{M}}_{g,n}(\mathbf{b}) = \left\{M\left|\begin{array}{l}
     M\ \mathrm{is\ a\ super \ hyperbolic\ surface\ of\ genus\ }g \mathrm{\ with} \\ \mathrm{NS\ geodesic\ boundary\ components\ of\ length}\ b_i \end{array}\right\}\right/\sim,
\]
where $\mathbf{b} = (b_1,\dots, b_n)$, and $M\sim M'$ if and only if there exists a (super) isometry that takes $M$ and its boundary $\partial M$ to $M'$ and $\partial M'$, respectively. Note that only boundaries of the NS type are included in the moduli space $\widehat{\mathcal{M}}_{g,n}(\mathbf{b})$. The entire moduli space, regardless of the spin types of boundary geodesics, consists of $2^{2g+n-1}$ connected components.

Similarly to the classical case, we define the super Teichm\"uller space $\widehat{\mathcal{T}}_{g,n}(\mathbf{b})$. This space is of dimension $(6g-6+2n|4g-4+2n)$, with coordinates given by the simple geodesic lengths $l_i$ and twist parameters $\tau_i$ for $i=1,\dots,3g-3+n$, and fermionic moduli $\psi_i$ and $\xi_i$ for $j=1,\dots,2g-2+n$\cite{Penner_2019}.

In generalizing the Kodaira-Spencer identification, the tangent space of the super Teichm\"uller space is identified as\cite{norbury2020enumerative}
\[
T_{M}\widehat{\mathcal{T}}_{g,n}(\mathbf{b}) = H^1(M_0,\mathcal{T}_{M})\oplus H^1(M_0,\mathcal{T}^{1/2}_{M}),
\]
where $M_0$ is the underlying (classical) hyperbolic surface in $M$. The first component corresponds to the tangent space of the classical Teichm\"uller space, while the second component is fermionic. The metric of the super Teichm\"uller space arising from this identification is given by\cite{stanford2019jt}
\[
    \bigwedge_{i=1}^{3g-3+n}(dl_i\wedge d\tau_i)\wedge \bigwedge_{j=1}^{2g-2+n}\left(\frac{1}{2}\exp\left(-\frac{l_{a_j}+l_{b_j}}{4}\right)\cosh\left(\frac{l_{c_j}}{4}\right)d\xi_j\wedge d\psi_j\right),
\]
where $a_j,b_j$ and $c_j$ are the indices for the three boundaries of the $j$-th pair of pants. The asymmetry between their lengths comes from Stanford and Witten's definition of the fermionic moduli $\xi_j$ and $\psi_j$. The volume form itself is indeed invariant under the exchange of these boundary component lengths.

The initial cases of $\widehat{V}_{g,n}(\mathbf{b})$, the volume of the super moduli space $\widehat{\mathcal{M}}_{g,n}(\mathbf{b})$ are given by
\begin{equation}\label{equ:2:2}
    \widehat{V}_{0,3}(\mathbf{b}) = 0, \widehat{V}_{1,1}(b) = \frac{1}{8}.
\end{equation} 
The vanishing of $V_{0,3}$ is due to that $\widehat{\mathcal{M}}_{0,3}$ has a greater fermionic dimension than its bosonic dimension, \cite{stanford2019jt}.

Similarly to the case of usual moduli space, a McShane identity holds for super hyperbolic surfaces:
\begin{theorem}[\cite{stanford2019jt}]
    Let $M$ be a super hyperbolic surface with genus $g$, a geodesic boundary component $\beta$ of length $b$, and other boundary components $\beta_i$ of length $b_i$ for $i=1,\dots,n-1$. Then,
    \[
    b = \sum_{(\alpha_j,\alpha_k)}\mathcal{D}(b,l(\alpha_j),l(\alpha_k)|\xi_{jk},\psi_{jk})+\sum_{i=1}^{n-1}\sum_{\alpha_j}\mathcal{R}(b,b_i,l(\alpha_j)|\xi_j,\psi_j),
    \]
    where the odd parameters $\xi_{jk}$, $\psi_{jk}$ are fermionic moduli for the pair of pants bounded by the thickening of simply closed geodesics $\beta$, $\alpha_j$ and $\alpha_k$, etc. Moreover, the function
    \[
        \mathcal{R}(x,y,z|\xi,\psi) = \mathcal{R}(x,y,z) - \xi\psi\frac{\exp\frac{x+z}{4}}{\cosh\frac{y}{4}}\widehat{\mathcal{R}}(x,y,z),
    \]
    where $\mathcal{R}(x,y,z)$ is the same as which defined by Mirzakhani, and
    \[
        \widehat{\mathcal{R}}(x,y,z) = \frac{1}{2}\widehat{H}(x+y,z)+\frac{1}{2}\widehat{H}(x-y,z),
    \]
    with
    \[
    \widehat{H}(x,y) = \frac{1}{2}\left(\left(\cosh\frac{x-y}{4}\right)^{-1}-\left(\cosh\frac{x+y}{4}\right)^{-1}\right).
    \]
    The function $\mathcal{D}(x,y,z|\xi,\psi)$ resembles a similar structure.
\end{theorem}
Integrating the identity over the two odd moduli eliminates the $\mathcal{R}(x,y,z)$ part, resulting in \eqref{equ:1:3} which involves only $\widehat{H}(x,y)$.

Analogously to Mirzakhani's case, the super volumes $\widehat{V}_{g,n}(\mathbf{b})$ are realized as integrals over the compact moduli space $\overline{\mathcal{M}}_{g,n}$, \cite{norbury2020enumerative}:
\begin{equation}\label{is:2}
\widehat{V}_{g,n}(\mathbf{b}) = \int_{\overline{\mathcal{M}}_{g,n}}\Theta_{g,n}\exp\left(2\pi^2\kappa_1 + \frac{1}{2}\sum b_i^2\psi_i\right),
\end{equation}
where the additional factor $\Theta_{g,n}\in H^{4g-4+2n}(\overline{\mathcal{M}}_{g,n})$ is called the \emph{Norbury class}. The Norbury class $\Theta_{g,n}$ can be reduced to the class $\Theta_{g,0}\in H^{4g-4}$, via the forgetful map $\pi:\overline{\mathcal{M}}_{g,n}\to \overline{\mathcal{M}}_{g,n-1}$. This reduction also explains the vanishing of $\widehat{V}_{0,n}(\mathbf{b})$ for $n\geq 3$.

%% file: Sec_3.tex
\section{Proof of Theorems \ref{thm:1:2} and \ref{thm:1:4}}
\label{Chp_2}
The proof of Theorems \ref{thm:1:2} and \ref{thm:1:4} requires us to consider the recursive relations \eqref{equ:1:1} and \eqref{equ:1:3} for a complex variable $b$. For $b'\in\mathbb{R}$, the functions $H(b',b)$ and $\widehat{H}(b',b)$ have poles on the lines $\mathrm{Im}(b) = \pm 2\pi$; however, it is safe to consider \eqref{equ:1:1} and \eqref{equ:1:3} for
\[
b\in \Omega:=\{z\in \mathbb{C}|-2\pi<\mathrm{Im}z<2\pi\}.
\]
By induction, equations \eqref{equ:1:1} and \eqref{equ:1:3} implies that the functions $V_{g,n}(b,\mathbf{b})$ and $\widehat{V}_{g,n}(b,\mathbf{b})$ are even and holomorphic in $x$ on the open region $\Omega$. We will now proceed by considering the two cases separately in the following subsections.
\subsection{The Recursion for Usual Moduli Space}
\label{Sec_2_1}
For $b,b_1,\dots, b_{n-1}\in \mathbb{R}$, we can consider $V_{g,n}(b\pm 2\pi \mathrm{i},\mathbf{b})$ as limits $\lim_{\beta\to \pm 2\pi_\mp}V_{g,n}(b + \beta \mathrm{i},\mathbf{b})$. Under this assumption, we prove Theorem \ref{thm:1:2} as follows.
\begin{lemma}\label{lem:3:1}
    The right-hand side of \eqref{equ:1:1} is bounded over the region $\{b = b_0+\beta \mathrm{i}|-B_0\leq b_0\leq B_0,0\leq \beta \leq 2\pi\}$ for any $B_0>0$ and any tuple $\mathbf{b}$. 
\end{lemma}
\begin{proof}
    We will demonstrate this for the first term on the right-hand side. We note that the only possible singularity in the integral
    \[
    \frac{1}{2}\int_0^\infty\int_0^\infty b'b''V_{g-1,n+1}(b',b'',\mathbf{b})H(b'+b'',b_0+\beta\mathrm{i})db'db''
    \]
    occurs when $b'+b'' = b_0$ and $\beta = 2\pi$. Therefore, if we separate the domain of the integral into two parts as follows:
    \[
    \begin{split}
        & \frac{1}{2}\iint_{\substack{b',b''\geq 0\\ b'+b''\geq 2B_0}} b'b''V_{g-1,n+1}(b',b'',\mathbf{b})H(b'+b'',b_0+\beta\mathrm{i})db'db''\\
        & + \frac{1}{2}\iint_{\substack{b',b''\geq 0\\ b'+b''\leq 2B_0}} b'b''V_{g-1,n+1}(b',b'',\mathbf{b})H(b'+b'',b_0+\beta\mathrm{i})db'db'',
    \end{split}
    \]
    then the first part of the integral is continuous and bounded for $b = b_0+\beta\mathrm{i}$ in this region.

    The second part is continuous for $\beta<2\pi$; when taking the limit $\beta\to 2\pi_-$, it converges to the sum of a principal value integral and a residue at the singularity $b'+b'' = b_0$:
    \[
    \begin{split}
        & \frac{1}{2}\iint_{\substack{b',b''\geq 0\\ b'+b''\leq 2B_0}} b'b''V_{g-1,n+1}(b',b'',\mathbf{b})H(b'+b'',b_0+\beta\mathrm{i})db'db'' \\
        & \overset{x = b'-b'',\ y = b'+b''}{=\joinrel=\joinrel=\joinrel=\joinrel=\joinrel=\joinrel=\joinrel=\joinrel=} \frac{1}{4}\int_0^{2 B_0}\int_{-y}^y\frac{y+x}{2}\frac{y-x}{2}V_{g-1,n+1}(\frac{y+x}{2},\frac{y-x}{2},\mathbf{b})H(y,b_0+\beta\mathrm{i})dxdy \\
        & \overset{\beta\to 2\pi}{\longrightarrow} \frac{1}{4}\mathrm{p.v.}\int_0^{2 B_0}\left(\int_{-y}^y\frac{y+x}{2}\frac{y-x}{2}V_{g-1,n+1}(\frac{y+x}{2},\frac{y-x}{2},\mathbf{b}) dx\right)H(y,b_0+2\pi \mathrm{i})dy \\
        & + \pi \mathrm{i}\mathrm{Res}_{y = b_0}\left(\frac{1}{4}\int_{-y}^y\frac{y+x}{2}\frac{y-x}{2}V_{g-1,n+1}(\frac{y+x}{2},\frac{y-x}{2},\mathbf{b}) dx\right)\frac{1}{1-\exp (y-b_0)/2},
    \end{split}
    \]
    and the residue is calculated as
    \[
        \frac{1}{2}\pi \mathrm{i}\int_{-b_0}^{b_0}\frac{b_0+x}{2}\frac{b_0-x}{2}V_{g-1,n+1}(\frac{b_0+x}{2},\frac{b_0-x}{2},\mathbf{b})dx.
    \]
    This implies that the second part is also bounded over the region $\{b = b_0+\beta \mathrm{i}|-B_0\leq b_0\leq B_0,0\leq \beta \leq 2\pi\}$ for any $B_0>0$ and any tuple $\mathbf{b}$. 
\end{proof}
\begin{proof}[Proof of Theorem \ref{thm:1:2}]
By Lemma \ref{lem:3:1}, the limit
\[
\lim_{\beta\to 2\pi_-}(b_0+\beta \mathrm{i})V_{g,n}(b_0+\beta \mathrm{i},\mathbf{b})
\]
exists for any $b_0\in\mathbb{R}$ and any tuple $\mathbf{b}$; the same holds for the limit as $\beta\to -2\pi_+$.

Moreover, the right-hand side of \eqref{equ:1:1} consists of integrals of functions that vanish exponentially and have only simple poles when $\mathrm{Im}(b) = \pm 2\pi\mathrm{i}$. Thus we can swap the order of integrals when integrating \eqref{equ:1:1} with respect to $b$ from $b_0 - 2\pi \mathrm{i}$ to $b_0+2\pi\mathrm{i}$. The left-hand side of the integration is
\[
    \int_{b_0-2\pi \mathrm{i}}^{b_0+2\pi \mathrm{i}}\frac{\partial}{\partial b}\left(bV_{g,n}(b,\mathbf{b})\right)db = (b_0+2\pi \mathrm{i})V_{g,n}(b_0+2\pi \mathrm{i},\mathbf{b}) - (b_0-2\pi \mathrm{i})V_{g,n}(b_0-2\pi \mathrm{i},\mathbf{b}),
\]
and the first term on the right-hand side of the integration is
\[
\begin{split}
    & \frac{1}{2}\int_{b_0-2\pi \mathrm{i}}^{b_0+2\pi \mathrm{i}}\int_0^\infty\int_0^\infty b'b''V_{g-1,n+1}(b',b'',\mathbf{b})H(b'+b'',b)db'db''db \\
    & = \frac{1}{2}\int_0^\infty\int_0^\infty b'b''V_{g-1,n+1}(b',b'',\mathbf{b})\left(\int_{b_0-2\pi \mathrm{i}}^{b_0+2\pi \mathrm{i}}H(b'+b'',b)db\right)db'db'' \\
    & = \frac{1}{2}\int_0^\infty\int_0^\infty b'b''V_{g-1,n+1}(b',b'',\mathbf{b})h(b'+b'',b_0)db'db'',
\end{split}
\]
where
\[
h(x,b): = \int_{b-2\pi \mathrm{i}}^{b+2\pi \mathrm{i}}H(x,b)db.
\]
Similarly, the other two terms on the right-hand side of the integration are
\[
\begin{split}
    & \frac{1}{2}\sum_{\substack{g_1+g_2 = g\\ \mathbf{b}_1\sqcup \mathbf{b}_2 = \mathbf{b}}}\int_0^\infty\int_0^\infty b'b''V_{g_1,n_1}(b',\mathbf{b}_1)V_{g_2,n_2}(b'',\mathbf{b}_2)h(b'+b'',b_0)db'db''\\
    & + \frac{1}{2}\sum_{j=1}^{n-1}\int_0^\infty b'V_{g,n-1}(b',\mathbf{b}\backslash b_j)\left(h(b',b_0+b_j)+h(b',b_0-b_j)\right)db'.
\end{split}
\]
The function $h(x,b)$ is computed as follows:
\[
\begin{split}
    & h(x,b) = \int_{b-2\pi\mathrm{i}}^{b+2\pi\mathrm{i}}H(x,b)db \overset{z = \mathrm{e}^{b/2}}{=\joinrel=} \int_{C_b} \left(\frac{1}{1+\mathrm{e}^{x/2}z}+\frac{z}{z+\mathrm{e}^{x/2}}\right)\frac{2dz}{z}\\
    = & 2\int_{C_b}\left(\frac{1}{z+\mathrm{e}^{x/2}}-\frac{1}{z+\mathrm{e}^{-x/2}}+\frac{1}{z}\right)dz = 
    \begin{cases}
        4\pi\mathrm{i}, & |x|<|b|,\\
        0, & x>|b|,\\
        8\pi\mathrm{i}, & x<-|b|,
    \end{cases}
\end{split}
\]
where $C_b$ is the circle $|z| = e^{b/2}$. That is to say, $h(x,b) = 4\pi\mathrm{i}\theta(|b|-x)$, where $\theta$ represents the Heaviside function. Substituting the expression for $h(x,b)$ into the integration of \eqref{equ:1:1}, we obtain
\[
    \begin{split}
        & \frac{1}{4\pi i}((b+2\pi i)V_{g,n}(b+2\pi i,\mathbf{b}) - (b-2\pi i)V_{g,n}(b-2\pi i,\mathbf{b}))\\
        & = \frac{1}{2}\iint_{\substack{b',b''\geq 0\\ b'+b''\leq |b|}} b'b''V_{g-1,n+1}(b',b'',\mathbf{b})db'db''\\
        & + \frac{1}{2}\sum_{\substack{g_1+g_2 = g\\ \mathbf{b}_1\sqcup \mathbf{b}_2 = \mathbf{b}}}\iint_{\substack{b',b''\geq 0\\ b'+b''\leq |b|}} b'b''V_{g_1,n_1}(b',\mathbf{b}_1)V_{g_2,n_2}(b'',\mathbf{b}_2)db'db''\\
        & + \frac{1}{2}\sum_{j=1}^{n-1}\left(\int_0^{|b+b_j|}+\int_0^{|b-b_j|}\right) b'V_{g,n-1}(b',\mathbf{b}\backslash b_j)db'.
    \end{split}
\]
Since the functions $V_{g,n}$ are even, we can eliminate all the absolute value symbols, resulting in \eqref{equ:1:2} for $b,b_1,\dots, b_{n-1}\in\mathbb{R}$.
\end{proof}
It is worth noting that our proof does not directly use the polynomiality of $V_{g,n}$. Instead, we can demonstrate the polynomiality of $V_{g,n}$ by utilizing \eqref{equ:1:2} along with the polynomial growth condition, which is guaranteed by Equation \eqref{equ:1:1}: 
\begin{Col}
    Suppose that functions $V_{g,n}(b,\mathbf{b})$ are holomorphic for $b\in \Omega$, exhibit polynomial growth in $b$, and satisfy the recursion relation \eqref{equ:1:2}. Then these functions are polynomials in all variables.
\end{Col}
\begin{proof}
    We will prove this by induction by $2g+n$. The base case $2g+n = 3$ holds following equation \eqref{equ:2:1}.

    Supposing that functions $V_{g',n'}$ are polynomials for any $2g'+n'<2g+n$, we aim to prove that $V_{g,n}$ is also a polynomial. By the induction assumption, the right-hand side of \eqref{equ:1:2} consists of integrals of polynomials over polygonal regions, which implies it is a polynomial. Therefore, there exists a polynomial function $\bar{V}_{g,n}$ such that
    \[
    \resizebox{\textwidth}{!}{$
    \begin{aligned}
    (b+2\pi \mathrm{i})V_{g,n}(b+2\pi \mathrm{i},\mathbf{b}) - (b-2\pi \mathrm{i})V_{g,n}(b-2\pi \mathrm{i},\mathbf{b}) = (b+2\pi \mathrm{i})\bar{V}_{g,n}(b+2\pi \mathrm{i},\mathbf{b}) - (b-2\pi \mathrm{i})\bar{V}_{g,n}(b-2\pi \mathrm{i},\mathbf{b}).
    \end{aligned}
    $}
    \]
    Denote $\tilde{V}(b;\mathbf{b}): = b(\bar{V}_{g,n}(b,\mathbf{b}) - V_{g,n}(b,\mathbf{b}))$. Then $\tilde{V}$ is holomorphic on $\Omega$, exhibit at most polynomial growth, and satisfies
    \[
    \tilde{V}(b+2\pi\mathrm{i}) = \tilde{V}(b-2\pi\mathrm{i}),
    \]
    for $b\in\mathbb{R}$. Consequently, $\tilde{V}$ extends to a holomorphic function on $\mathbb{C}$, with $\tilde{V}(0) = 0$, exhibits polynomial growth and period $4\pi\mathrm{i}$. It follows that $\tilde{V}\equiv 0$. Therefore, $V_{g,n} = \tilde{V}_{g,n}$ is a polynomial on $b$. By symmetry, $V_{g,n}$ is also a polynomial in the other variables. 
\end{proof}
As a function over $\mathbb{C}$, the difference between the two sides of \eqref{equ:1:2} is a polynomial on $b$ and vanishes when restricting to $b\in\mathbb{R}$. Therefore, we conclude:
\begin{Col}
    Equation \eqref{equ:1:2} holds for any $b$, $b_1,\dots,b_{n-1}\in\mathbb{C}$.
\end{Col}
\subsection{The Recursion for Super Moduli Space}
\label{Sec_2_2}
Similarly to the case of usual moduli spaces, we consider $\widehat{V}_{g,n}(b\pm 2\pi\mathrm{i},\mathbf{b})$ as limits. Similar to Lemma \ref{lem:3:1}, for $b = b_0+\beta\mathrm{i}$, the limit of the right-hand side of \eqref{equ:1:3} is bounded as $\beta \to \pm 2\pi$ for $b_0$ in a bounded region. For $b,b_1,\dots,b_{n-1}\in\mathbb{R}$, Theorem \ref{thm:1:4} is proved as follows.
\begin{proof}[Proof of Theorem \ref{thm:1:4}]
    By Theorem \ref{Thm:1:3}, the limits
    \[
    \lim_{\beta\to 2\pi_-}(b_0+\beta\mathrm{i})\widehat{V}(b_0+\beta\mathrm{i},\mathbf{b})\ \text{and}\ \lim_{\beta\to -2\pi_+}(b_0+\beta\mathrm{i})\widehat{V}(b_0+\beta\mathrm{i},\mathbf{b})
    \]
    exist for any real $b_0,b_1,\dots, b_{n-1}$. Evaluating equation \eqref{equ:1:3} at $b = b_0 + 2\pi\mathrm{i}$ and $b = b_0 - 2\pi\mathrm{i}$, respectively, and taking the sum of these, the left-hand side of the new equation is
    \[
    2\pi\left((b+2\pi i)\widehat{V}_{g,n}(b+2\pi i,\mathbf{b}) + (b-2\pi i)\widehat{V}_{g,n}(b-2\pi i,\mathbf{b})\right).
    \]
    To evaluate the right-hand side of the equation, note that
    \[
    \cosh(z\pm \frac{\pi\mathrm{i}}{2}) = \pm \mathrm{i}\sinh(z),
    \]
    for any $\mathrm{z}\in\mathbb{C}$. Utilizing this fact, the first term on the right-hand side can be expressed as
    \[
    \resizebox{\textwidth}{!}{$
    \begin{aligned}
        & \frac{1}{2}\lim_{\epsilon\to 0_+}\int_0^\infty\int_0^\infty b'b''\widehat{V}_{g-1,n+1}(b',b'',\mathbf{b})(\widehat{H}(b'+b'',b_0+(2\pi-\epsilon)\mathrm{i})+\widehat{H}(b'+b'',b_0-(2\pi-\epsilon)\mathrm{i}))db'db'' \\
        & = \frac{\mathrm{i}}{4}\lim_{\epsilon\to 0_+}\int_0^\infty\int_0^\infty b'b''\widehat{V}_{g-1,n+1}(b',b'',\mathbf{b})\cdot\left(\sinh\left(\frac{b'+b''+b_0+\epsilon\mathrm{i}}{4}\right)^{-1}\right.\\
        & \left. - \sinh\left(\frac{b'+b''-b_0+\epsilon\mathrm{i}}{4}\right)^{-1} - \sinh\left(\frac{b'+b''+b_0-\epsilon\mathrm{i}}{4}\right)^{-1} + \sinh\left(\frac{b'+b''-b_0-\epsilon\mathrm{i}}{4}\right)^{-1}\right)db'db''\\
        & = \frac{\mathrm{i}}{4}\lim_{\epsilon\to 0_+}\left(\int_{-\epsilon\mathrm{i}}^{\infty-\epsilon\mathrm{i}} - \int_{\epsilon\mathrm{i}}^{\infty+\epsilon\mathrm{i}}\right)\int_0^\infty b'b''\widehat{V}_{g-1,n+1}(b',b'',\mathbf{b})\\
        & \cdot\left(\sinh\left(\frac{b'+b''-b_0}{4}\right)^{-1}-\sinh\left(\frac{b'+b''+b_0}{4}\right)^{-1}\right)db'db'' \\
        & = \frac{\mathrm{i}}{4}\lim_{\substack{\epsilon\to 0_+\\N\to\infty}}\int_0^\infty\int_{C_{\epsilon,N}} b'b''\widehat{V}_{g-1,n+1}(b',b'',\mathbf{b})\left(\sinh\left(\frac{b'+b''-b_0}{4}\right)^{-1}-\sinh\left(\frac{b'+b''+b_0}{4}\right)^{-1}\right)db''db',
    \end{aligned}
    $}
    \]
    where $C_{\epsilon,N}$ denotes the rectangular path given by
    \[
    -\epsilon\mathrm{i}\to N-\epsilon\mathrm{i}\to N+\epsilon\mathrm{i}\to \epsilon\mathrm{i}\to -\epsilon\mathrm{i};
    \]
    see Figure \ref{fig:1} below.
    \begin{figure}[H]
    \centering
    \begin{tikzpicture}
        \draw[thick,->] (0,-0.5) -- (3.5,-0.5);
        \draw[thick,-] (3.5,-0.5) -- (7,-0.5) node[anchor=west]{$N-\epsilon\mathrm{i}$};
        \draw[thick,->] (7,-0.5) -- (7,0);
        \draw[thick,-] (7,0) -- (7,0.5) node[anchor=west]{$N+\epsilon\mathrm{i}$};
        \draw[thick,->] (7,0.5) -- (3.5,0.5);
        \draw[thick,-] (3.5,0.5) -- (0,0.5) node[anchor=east]{$\epsilon\mathrm{i}$};
        \draw[thick,->] (0,0.5) -- (0,0);
        \draw[thick,-] (0,0) -- (0,-0.5) node[anchor=east]{$-\epsilon\mathrm{i}$};
    \end{tikzpicture}
    \caption{The path $C_{\epsilon,N}$.}
    \label{fig:1}
    \end{figure}
    As the parameter $\epsilon$ approaches $0$ and $N$ becomes sufficiently large, the function inside the integral above possesses a unique pole at $b'' = |b_0|-b'$ if $b'<|b_0|$, or no poles if $b'>|b_0|$. We assume $b_0\geq 0$, as the other case can be verified similarly. We simplify the integral above by evaluating the residue at this simple pole:
    \[
    \resizebox{\textwidth}{!}{$
    \begin{aligned}
        & \frac{1}{2}\lim_{\epsilon\to 0_+}\int_0^\infty\int_0^\infty b'b''\widehat{V}_{g-1,n+1}(b',b'',\mathbf{b})(\widehat{H}(b'+b'',b_0+(2\pi-\epsilon)\mathrm{i})+\widehat{H}(b'+b'',b_0-(2\pi-\epsilon)\mathrm{i}))db'db''\\
        & = \frac{\mathrm{i}}{4}\int_0^{b_0} 2\pi\mathrm{i}\underset{b''= b_0-b'}{\mathrm{Res}}b'b''\widehat{V}_{g-1,n+1}(b',b'',\mathbf{b})\left(\sinh\left(\frac{b'+b''-b_0}{4}\right)^{-1}-\sinh\left(\frac{b'+b''+b_0}{4}\right)^{-1}\right)db'\\
        & = -2\pi\int_0^{b_0} b'(b_0-b')\widehat{V}_{g-1,n+1}(b',b_0-b',\mathbf{b})db'.
    \end{aligned}
    $}
    \]
    Similarly, the other two terms on the right-hand side simplify to
    \[
    -2\pi\sum_{\substack{g_1+g_2 = g\\ \mathbf{b}_1\sqcup \mathbf{b}_2 = \mathbf{b}}}\int_0^{b_0} b'(b_0-b')\widehat{V}_{g_1,n_1}(b',\mathbf{b}_1)\widehat{V}_{g_2,n_2}(b_0-b',\mathbf{b}_2)db'
    \]
    and
    \[
    -2\pi\sum_{j=1}^{n-1}\left((b_0+b_j)\widehat{V}_{g,n-1}(b_0+b_j,\mathbf{b}\backslash b_j)+(b_0-b_j)\widehat{V}_{g,n-1}(b_0-b_j,\mathbf{b}\backslash b_j)\right),
    \]
    respectively. Dividing these expressions by $-4\pi$ yields equation \eqref{equ:1:4}.
\end{proof}
Let us show the polynomiality of $\widehat{V}_{g,n}$ utilizing \eqref{equ:1:4} along with the polynomial growth condition guaranteed by \eqref{equ:1:3}.
\begin{Col}
    Suppose that functions $\widehat{V}_{g,n}(b,\mathbf{b})$ are holomorphic for $b\in \Omega$, exhibit polynomial growth on $b$, and satisfy the recursion relation \eqref{equ:1:4}. Then these functions are polynomials in all variables.
\end{Col}
\begin{proof}
    We will show this using induction on $2g+n$. The base case $2g+n = 3$ is guaranteed by equation \eqref{equ:2:2}.

    Assuming that functions $\widehat{V}_{g',n'}$ are polynomials for any $2g'+n'<2g+n$, we aim to prove that $\widehat{V}_{g,n}$ is also a polynomial. By the induction assumption, the right-hand side of \eqref{equ:1:4} is a polynomial. Thus, there exists a polynomial function $\bar{\widehat{V}}_{g,n}$ such that
    \[
    \resizebox{\textwidth}{!}{$
    \begin{aligned}
    (b+2\pi \mathrm{i})\widehat{V}_{g,n}(b+2\pi \mathrm{i},\mathbf{b}) + (b-2\pi \mathrm{i})\widehat{V}_{g,n}(b-2\pi \mathrm{i},\mathbf{b}) = (b+2\pi \mathrm{i})\bar{\widehat{V}}_{g,n}(b+2\pi \mathrm{i},\mathbf{b}) + (b-2\pi \mathrm{i})\bar{\widehat{V}}_{g,n}(b-2\pi \mathrm{i},\mathbf{b}).
    \end{aligned}
    $}
    \]
    Similarly, we denote $\tilde{\widehat{V}} = b(\bar{\widehat{V}}_{g,n} - \widehat{V}_{g,n})$. Then
    \[
    \tilde{\widehat{V}}(b+2\pi \mathrm{i}) = -\tilde{\widehat{V}}(b-2\pi \mathrm{i}),
    \]
    for $b\in\mathbb{R}$. The rest of the proof follows similarly, noting that $\tilde{\widehat{V}}$ can be extended to a holomorphic function of period $8\pi\mathrm{i}$. It follows that \( \tilde{\widehat{V}} \) vanishes. Consequently, \( \widehat{V}_{g,n} \) is a polynomial, which completes the proof.
\end{proof}
\begin{Col}
    Equation \eqref{equ:1:4} holds for any $b$, $b_1,\dots,b_{n-1}\in\mathbb{C}$.
\end{Col}

%% file: Sec_4.tex
\section{The Laplace form of original and new recursion formulas}
\label{Chp_3}
By applying a discrete Laplace Transform to the cut-and-join equation for Catalan numbers \cite{dumitrescu2015lectures} and Hurwitz numbers \cite{mulase2012spectral}, one established the corresponding \emph{topological recursion}. In this section, a similar strategy is applied to the recursion equations \eqref{equ:1:1} through \eqref{equ:1:4} for the moduli space volume functions. Due to the continuous nature of these recursion formulas, we consider a continuous Laplace transform for them.
\subsection{The Recursion for Usual Moduli Space}
\label{Sec_3_1}
Recall that the Laplace transform $F_{g,n}^V$ is defined as
\[
F_{g,n}^V(\mathbf{t}) = \int_{\mathbb{R}_+^n}V_{g,n}(\mathbf{b})\mathrm{e}^{-\mathbf{b}\cdot\mathbf{t}}d\mathbf{b}.
\]
For example,
\[
F_{0,3}^V(t_1,t_2,t_3) = \int_0^\infty\int_0^\infty\int_0^\infty \mathrm{e}^{-b_1t_1-b_2t_2-b_3t_3}db_1db_2db_3 = \frac{1}{t_1 t_2 t_3},
\]
and
\[
F_{1,1}^V(t) = \int_0^\infty \frac{(b^2+4\pi^2)\mathrm{e}^{-bt}}{48}db = \frac{1}{24t^3}+\frac{\pi^2}{12 t}.
\]
Note that the volume functions $V_{g,n}$ are even polynomials, as we have realized in Section \ref{Chp_2}. Thus, we have a first conclusion for their Laplace transforms:
\begin{prop}
    The Laplace transforms $F_{g,n}^V(\mathbf{t})$ of Weil-Petersson volumes are odd polynomials in $t_1^{-1},\dots, t_n^{-1}$.
\end{prop}
The proof of Theorem \ref{thm:1:5} requires the following lemmas, which are elementary and will be proved in Appendix \ref{Append:1}.
\begin{lemma}\label{lem:4:1}
    Let $b,t\in\mathbb{C}_{+}$, i.e., $\mathrm{Re}(b)$, $\mathrm{Re}(t)>0$. Then the following integrals and series converge and are equal to each other:
    \begin{equation}\label{equ:3:1}
    \int_0^\infty \frac{\mathrm{e}^{-tx}dx}{1+\mathrm{e}^{(b+x)/2}} = -\sum_{k=1}^\infty \frac{(-1)^k\mathrm{e}^{-bk/2}}{t+k/2}.
    \end{equation}
    Additionally, if $t\notin \frac{1}{2}\mathbb{Z}$, then the following integrals and series also converge and are equal to each other:
    \begin{equation}\label{equ:3:2}
    \int_0^\infty \frac{\mathrm{e}^{-tx}dx}{1+\mathrm{e}^{(b-x)/2}} = \frac{2\pi \mathrm{e}^{-tb}}{\sin 2\pi t}-\sum_{k=1}^\infty \frac{(-1)^k\mathrm{e}^{-bk/2}}{t-k/2}.
    \end{equation}
\end{lemma}
\begin{lemma}\label{lem:4:2}
    Let $b,t_1,t_2\in \mathbb{C}_{+}$. Then the following integrals and series converge and are equal to each other:
    \begin{equation}\label{equ:3:3}
    \int_0^\infty\int_0^\infty\frac{\mathrm{e}^{-t_1x_1-t_2x_2}dx_1dx_2}{1+\mathrm{e}^{(b+x_1+x_2)/2}} = -\sum_{k=1}^\infty\frac{(-1)^k\mathrm{e}^{-bk/2}}{(t_1+k/2)(t_2+k/2)},
    \end{equation}
    \begin{equation}\label{equ:3:4}
    \int_0^\infty\int_0^\infty\frac{\mathrm{e}^{-t_1x_1-t_2x_2}dx_1dx_2}{1+\mathrm{e}^{(b-x_1+x_2)/2}} = \frac{2\pi \mathrm{e}^{-t_1b}}{(t_1+t_2)\sin 2\pi t_1} -\sum_{k=1}^\infty\frac{(-1)^k\mathrm{e}^{-bk/2}}{(t_1-k/2)(t_2+k/2)},
    \end{equation}
    \begin{equation}\label{equ:3:5}
    \int_0^\infty\int_0^\infty\frac{\mathrm{e}^{-t_1x_1-t_2x_2}dx_1dx_2}{1+\mathrm{e}^{(b+x_1-x_2)/2}} = \frac{2\pi \mathrm{e}^{-t_2b}}{(t_1+t_2)\sin 2\pi t_2} -\sum_{k=1}^\infty\frac{(-1)^k\mathrm{e}^{-bk/2}}{(t_1+k/2)(t_2-k/2)},
    \end{equation}
    \begin{equation}\label{equ:3:6}
    \resizebox{0.92\textwidth}{!}{$
    \begin{aligned}
    \int_0^\infty\int_0^\infty\frac{\mathrm{e}^{-t_1x_1-t_2x_2}dx_1dx_2}{1+\mathrm{e}^{(b-x_1-x_2)/2}} = \frac{2\pi}{t_1-t_2}\left(\frac{\mathrm{e}^{-t_2b}}{\sin 2\pi t_2} - \frac{\mathrm{e}^{-t_1b}}{\sin 2\pi t_1}\right)-\sum_{k=1}^\infty\frac{(-1)^k\mathrm{e}^{-bk/2}}{(t_1-k/2)(t_2-k/2)}.
    \end{aligned}
    $}
    \end{equation}
    Here we assume that $t_1\notin \frac{1}{2}\mathbb{Z}$ for equations \eqref{equ:3:4} and \eqref{equ:3:6}, $t_2\notin \frac{1}{2}\mathbb{Z}$ for equations \eqref{equ:3:5} and \eqref{equ:3:6}, and $t_1\neq t_2$ for equation \eqref{equ:3:6}.
\end{lemma}
\begin{lemma}\label{lem:4:3}
    Let $\mathsf{H}_t$ denote taking the holomorphic part of a Laurent series at $t = 0$. Then for any even number $m\in \mathbb{N}$ and $t\neq \frac{1}{2}\mathbb{Z}$, the following series converge and are equal to each other:
        \begin{equation}\label{equ:4:3}
        \sum_{k=1}^\infty\frac{(-1)^k2t(2/k)^m}{t^2-k^2/4} = \mathsf{H}_t\left(\frac{2\pi}{t^m\sin 2\pi t}\right).
        \end{equation}
\end{lemma}
\begin{lemma}\label{lem:4:4}
    For any odd number $m\in\mathbb{N}$, the principal part
    \[
    \mathsf{P}_t\left(\frac{st^{-m}}{s^2-t^2}\right) = \frac{st^{-m} - ts^{-m}}{s^2-t^2}.
    \]
    For any even number $m\in\mathbb{N}$, the principal part
    \[
    \mathsf{P}_t\left(\frac{st^{-m}}{s^2-t^2}\right) = \frac{st^{-m} - s^{1-m}}{s^2-t^2}.
    \]
\end{lemma}
We return to the proof of Theorem \ref{thm:1:5}.
\begin{proof}[Proof of Theorem \ref{thm:1:5}, statement (1)]
    We will show that Equation \eqref{equ:1:5} can be obtained by multiplying the cut-and join equation \eqref{equ:1:1} with $\mathrm{e}^{-bt -\mathbf{b}\cdot\mathbf{t}}$ and integrating the result by $b$, $b_1$,\dots, $b_{n-1}$ from $0$ to $\infty$ respectively. Since $V_{g,n}(b,\mathbf{b})$ is a polynomial, these integrals absolutely converge with the assumption $t$, $t_1,\dots, t_{n-1}\in \mathbb{C}_+$. The left-hand side of the resulting equation reads as
    \[
    \begin{split}
        & \int_{\mathbb{R}_+^{n-1}}\int_0^\infty \frac{\partial}{\partial b}(b V_{g,n}(b,\mathbf{b}))\mathrm{e}^{-bt -\mathbf{b}\cdot\mathbf{t}} dbd\mathbf{b}\\
        & = -\int_{\mathbb{R}_+^{n-1}}\int_0^\infty b V_{g,n}(b,\mathbf{b})\frac{\partial}{\partial b}(\mathrm{e}^{-bt -\mathbf{b}\cdot\mathbf{t}}) dbd\mathbf{b}\\
        & = t\int_{\mathbb{R}_+^{n-1}}\int_0^\infty b V_{g,n}(b,\mathbf{b})\mathrm{e}^{-bt -\mathbf{b}\cdot\mathbf{t}} dbd\mathbf{b}\\
        & = -t\frac{\partial}{\partial t}\int_{\mathbb{R}_+^{n-1}}\int_0^\infty V_{g,n}(b,\mathbf{b})\mathrm{e}^{-bt -\mathbf{b}\cdot\mathbf{t}} dbd\mathbf{b} = -t\frac{\partial}{\partial t}F_{g,n}^V(t,\mathbf{t}).
    \end{split}
    \]
    The first term on the right-hand side equals to
    \[
    \begin{split}
        & \frac{1}{2}\int_{\mathbb{R}_+^{n-1}}\int_0^\infty\int_0^\infty\int_0^\infty b'b''V_{g-1,n+1}(b',b'',\mathbf{b})H(b'+b'',b)\mathrm{e}^{-bt -\mathbf{b}\cdot\mathbf{t}}db'db''dbd\mathbf{b} \\
        & = \frac{1}{2}\int_{\mathbb{R}_+^{n-1}}\int_0^\infty\int_0^\infty\left(\int_0^\infty H(b'+b'',b)\mathrm{e}^{-bt}db\right)b'b''V_{g-1,n+1}(b',b'',\mathbf{b})\mathrm{e}^{-\mathbf{b}\cdot\mathbf{t}}db'db''d\mathbf{b}.
    \end{split}
    \]
    Applying Lemma \ref{lem:4:1}, the integral over $b$ transforms to
    \[
    \begin{split}
        & \int_0^\infty H(b'+b'',b)\mathrm{e}^{-bt}db\\
        & = \int_0^\infty \frac{\mathrm{e}^{-bt}}{1+\mathrm{e}^{(b'+b''+b)/2}}+\frac{\mathrm{e}^{-bt}}{1+\mathrm{e}^{(b'+b''-b)/2}}db \\
        & = \frac{2\pi\mathrm{e}^{-t(b'+b'')}}{\sin 2\pi t} - \sum_{k=1}^\infty \frac{(-1)^k\mathrm{e}^{-(b'+b'')k/2}}{t+k/2} - \sum_{k=1}^\infty \frac{(-1)^k\mathrm{e}^{-(b'+b'')k/2}}{t-k/2} \\
        & = \frac{2\pi\mathrm{e}^{-t(b'+b'')}}{\sin 2\pi t} - \sum_{k=1}^\infty \frac{(-1)^k2t\mathrm{e}^{-(b'+b'')k/2}}{t^2- k^2/4}.
    \end{split}
    \]
    Therefore, the first term on the right-hand side breaks into two parts; the first part is simplified as
    \[
    \begin{split}
        & \frac{1}{2}\int_{\mathbb{R}_+^{n-1}}\int_0^\infty\int_0^\infty\frac{2\pi\mathrm{e}^{-t(b'+b'')}}{\sin 2\pi t}b'b''V_{g-1,n+1}(b',b'',\mathbf{b})\mathrm{e}^{-\mathbf{b}\cdot\mathbf{t}}db'db''d\mathbf{b}\\
        & = \frac{\pi}{\sin 2\pi t}\int_{\mathbb{R}_+^{n-1}}\int_0^\infty\int_0^\infty b'b''V_{g-1,n+1}(b',b'',\mathbf{b})\mathrm{e}^{-(b'+b'')t}\mathrm{e}^{-\mathbf{b}\cdot\mathbf{t}}db'db''d\mathbf{b} \\
        & = \frac{\pi}{\sin 2\pi t}\left(\left.\frac{\partial^2}{\partial t\partial t'}F_{g-1,n+1}^V(t,t',\mathbf{t})\right)\right|_{t' = t},
    \end{split}
    \]
    while the second part is simplified as
    \[
    \begin{split}
        & -\frac{1}{2}\int_{\mathbb{R}_+^{n-1}}\int_0^\infty\int_0^\infty\sum_{k=1}^\infty \frac{(-1)^k2t\mathrm{e}^{-(b'+b'')k/2}}{t^2- k^2/4}b'b''V_{g-1,n+1}(b',b'',\mathbf{b})\mathrm{e}^{-\mathbf{b}\cdot\mathbf{t}}db'db''d\mathbf{b}\\
        & = -\sum_{k=1}^\infty \frac{(-1)^kt}{t^2- k^2/4}\int_{\mathbb{R}_+^{n-1}}\int_0^\infty\int_0^\infty b'b''V_{g-1,n+1}(b',b'',\mathbf{b})\mathrm{e}^{-(b'+b'')k/2}\mathrm{e}^{-\mathbf{b}\cdot\mathbf{t}}db'db''d\mathbf{b} \\
        & = -\sum_{k=1}^\infty \frac{(-1)^kt}{t^2- k^2/4}\left(\left.\frac{\partial^2}{\partial t\partial t'}F_{g-1,n+1}^V(t,t',\mathbf{t})\right)\right|_{(t,t') = (k/2,k/2)}\\
        & = -\mathsf{H}_t\left(\frac{\pi}{\sin 2\pi t}\left(\left.\frac{\partial^2}{\partial t\partial t'}F_{g-1,n+1}^V(t,t',\mathbf{t})\right)\right|_{t' = t}\right).
    \end{split}
    \]
    Here $F^V_{g-1,n+1}(t,t',\mathbf{t})$ is an odd polynomial on both $1/t$ and $1/t'$, thus Lemma \ref{lem:4:3} is applied to the simplifying process. Combining these two parts, we simplify the first term on the right-hand side as
    \[
        \mathsf{P}_t\left(\frac{\pi}{\sin 2\pi t}\left(\left.\frac{\partial^2}{\partial t\partial t'}F_{g-1,n+1}^V(t,t',\mathbf{t})\right)\right|_{t' = t}\right).
    \]
    Similarly, the second term on the right-hand side is simplified as
    \[
        \sum_{\substack{g_1+g_2 = g\\\mathbf{t}_1\sqcup \mathbf{t}_2 = \mathbf{t}}}\mathsf{P}_t\left(\frac{\pi}{\sin 2\pi t}\left(\frac{\partial}{\partial t}F_{g_1,n_1}^V(t,\mathbf{t}_1)\right)\left(\frac{\partial}{\partial t}F_{g_2,n_2}^V(t,\mathbf{t}_2)\right)\right).
    \]
    For the third term, we apply Lemma \ref{lem:4:2} instead and transform the term as follows:
    \[
    \resizebox{\textwidth}{!}{$
    \begin{aligned}
        & \frac{1}{2}\sum_{j=1}^{n-1}\int_{\mathbb{R}_+^{n-1}}\int_0^\infty\int_0^\infty b'V_{g,n-1}(b',\mathbf{b}\backslash b_j)(H(b',b+b_j)+H(b',b-b_j))\mathrm{e}^{-bt -\mathbf{b}\cdot\mathbf{t}} dbdb'd\mathbf{b} \\
        & = \frac{1}{2}\sum_{j=1}^{n-1}\int_{\mathbb{R}_+^{n-2}}\int_0^\infty\left(\int_0^\infty\int_0^\infty (H(b',b+b_j)+H(b',b-b_j))\mathrm{e}^{-bt- b_jt_j}dbdb_j\right)b'V_{g,n-1}(b',\mathbf{b}\backslash b_j)\mathrm{e}^{-\mathbf{b}\cdot\mathbf{t}+b_jt_j}db'd(\mathbf{b}\backslash b_j) \\
        & = \frac{1}{2}\sum_{j=1}^{n-1}\int_{\mathbb{R}_+^{n-2}}\int_0^\infty\left(\frac{4\pi t_j\mathrm{e}^{-b't}}{(t_j^2-t^2)\sin 2\pi t} + \frac{4\pi t\mathrm{e}^{-b't_j}}{(t^2-t_j^2)\sin 2\pi t_j}\right) V_{g,n-1}(b',\mathbf{b}\backslash b_j)b' \mathrm{e}^{-\mathbf{b}\cdot \mathbf{t}+b_jt_j}db'd(\mathbf{b}\backslash b_j) \\
        & - \frac{1}{2}\sum_{j=1}^{n-1}\int_{\mathbb{R}_+^{n-2}}\int_0^\infty\sum_{k=1}^\infty\left(\frac{4(-1)^ktt_j\mathrm{e}^{-b'k/2}}{(t_j^2-t^2)(t^2-k^2/4)} + \frac{4(-1)^ktt_j\mathrm{e}^{-b'k/2}}{(t^2-t_j^2)(t_j^2-k^2/4)}\right) V_{g,n-1}(b',\mathbf{b}\backslash b_j)b' \mathrm{e}^{-\mathbf{b}\cdot \mathbf{t}+b_jt_j}db'd(\mathbf{b}\backslash b_j).
    \end{aligned}
    $}
    \]
    Note that the expression above consists of four parts, the first two parts are simplified as
    \[
    \sum_{j=1}^{n-1}\frac{2\pi t_j}{(t^2 - t_j^2)\sin 2\pi t}\frac{\partial}{\partial t}F_{g,n-1}^V(t,\mathbf{t}\backslash t_j)
    \]
    and
    \[
    -\sum_{j=1}^{n-1}\frac{2\pi t}{(t^2 - t_j^2)\sin 2\pi t_j}\frac{\partial}{\partial t_j}F_{g,n-1}^V(\mathbf{t}),
    \]
    respectively, while the last two parts are simplified as
    \[
    -\sum_{j=1}^{n-1}\frac{1}{t^2 - t_j^2}\mathsf{H}_t\left(\frac{2\pi t_j}{\sin 2\pi t}\frac{\partial}{\partial t}F_{g,n-1}^V(t,\mathbf{t}\backslash t_j)\right)
    \]
    and
    \[
    \sum_{j=1}^{n-1}\frac{1}{t^2 - t_j^2}\mathsf{H}_{t_j}\left(\frac{2\pi t}{\sin 2\pi t_j}\frac{\partial}{\partial t_j}F_{g,n-1}^V(\mathbf{t})\right)
    \]
    respectively, due to Lemma \ref{lem:4:3}. Collecting these parts, we simplify the third term on the right-hand side to
    \[
    \sum_{j=1}^{n-1}\frac{t_j\mathsf{P}_t\left(\frac{2\pi}{\sin 2\pi t}\frac{\partial}{\partial t}F_{g,n-1}^V(t,\mathbf{t}\backslash t_j)\right) - t\mathsf{P}_{t_j}\left(\frac{2\pi}{\sin 2\pi t_j}\frac{\partial}{\partial t_j}F_{g,n-1}^V(\mathbf{t})\right)}{t^2 - t_j^2}.
    \]
    By denoting
    \[
    P(t^{-1}): = \mathsf{P}_t\left(\frac{2\pi}{\sin 2\pi t}\frac{\partial}{\partial t}F_{g,n-1}^V(t,\mathbf{t}\backslash t_j)\right),
    \]
    a polynomial on $t^{-1}$, the third term can be expressed as
    \[
    \sum_{j=1}^{n-1}\frac{t_jP(t^{-1}) - tP(t_j^{-1})}{t^2 - t_j^2}.
    \]
    Therefore, by Lemma \ref{lem:4:4}, the third term on the right-hand side is finally simplified to
    \[
    \sum_{j=1}^{n-1}\mathsf{P}_t\frac{t_j\mathsf{P}_t\left(\frac{2\pi}{\sin 2\pi t}\frac{\partial}{\partial t}F_{g,n-1}^V(t,\mathbf{t}\backslash t_j)\right)}{t^2 - t_j^2} = \sum_{j=1}^{n-1}\mathsf{P}_t \left(\frac{2\pi t_j}{(t^2 - t_j^2)\sin 2\pi t}\frac{\partial}{\partial t}F_{g,n-1}^V(t,\mathbf{t}\backslash t_j)\right).
    \]
    These conclude equation \eqref{equ:1:5}.
\end{proof}

Eynard and Orantin, \cite{eynard2007weil}, considered a different form of the Laplace transform, related to $F_{g,n}^V$ as follows:
\[
    F_{g,n}^{EO}(\mathbf{t}) = \frac{(-1)^n\partial^n}{\partial t_1\dots\partial t_n}F_{g,n}^V(\mathbf{t}).
\]
For convenience, additional notation is introduced:
\[
    F_{0,2}^{EO}(t_1,t_2) = \frac{1}{(t_1-t_2)^2}.
\]
They derived a topological recursion for $F_{g,n}^{EO}$:
\begin{theorem}[\cite{eynard2007weil}]\label{thm:eo}
    The function $F_{g,n}^{EO}$ satisfies the following recursion formula:
    \[
    \resizebox{\textwidth}{!}{$
    \begin{aligned}
        F_{g,n}^{EO}(t,\mathbf{t})= \mathrm{Res}_s\frac{\pi}{(t^2-s^2)\sin 2\pi s}\left(F_{g-1,n+1}^{EO}(s,-s,\mathbf{t}) + \sum_{\substack{g_1+g_2 = g\\ \mathbf{t}_1\sqcup \mathbf{t}_2 = \mathbf{t}}}F_{g_1,n_1}^{EO}(s,\mathbf{t}_1)F_{g_2,n_2}^{EO}(-s,\mathbf{t}_2)\right).
    \end{aligned}
    $}
    \]
\end{theorem}
We give a proof of Theorem \ref{thm:eo} based on equation\eqref{equ:1:5}.
\begin{lemma}\label{lem:4:5}
    Suppose that
    \[
    f(z) = \sum_{k=-\infty}^\infty a_k z^{2k}
    \]
    is a Laurent series with only odd-order terms. Then for any $z$ in the converging region of its principal part, we have
    \[
    \mathrm{Res}_s \frac{f(s)}{z^2-s^2} = \frac{\mathsf{P}_z f(z)}{z}.
    \]
\end{lemma}
\begin{proof}[Proof of Theorem \ref{thm:eo}]
Note that $F_{g,n}^{EO}$ is an even function. According to Lemma \ref{lem:4:5}, it suffices to prove that
\begin{equation}\label{Equ_Proof_1}
    \begin{split}
        F_{g,n}^{EO}(t,\mathbf{t})= \frac{1}{t}\mathsf{P}_t\frac{\pi}{\sin 2\pi t}\left(F_{g-1,n+1}^{EO}(t,-t,\mathbf{t}) + \sum_{\substack{g_1+g_2 = g\\ \mathbf{t}_1\sqcup \mathbf{t}_2 = \mathbf{t}}}F_{g_1,n_1}^{EO}(t,\mathbf{t}_1)F_{g_2,n_2}^{EO}(-t,\mathbf{t}_2)\right).
    \end{split}
\end{equation}
To derive equation \eqref{Equ_Proof_1}, we take partial derivatives $\frac{(-1)^{n-1}\partial^{n-1}}{\partial t_1\dots\partial t_{n-1}}$ of both sides of \eqref{equ:1:5}. The left-hand side then becomes
\[
t\frac{(-1)^n\partial^n}{\partial t\partial t_1\dots \partial t_{n-1}}F_{g,n}^V(t,\mathbf{t}) = F_{g,n}^{EO}(t,\mathbf{t}).
\]
For the first term on the right-hand side, we have
\[
    \begin{split}
        & \frac{(-1)^{n-1}\partial^{n-1}}{\partial t_1\dots\partial t_{n-1}}\mathsf{P}_t\left(\frac{\pi}{\sin 2\pi t}\left.\frac{\partial^2}{\partial t\partial t'}F_{g-1,n+1}^V(t,t',\mathbf{t})\right|_{t' = t}\right) \\
        & = \mathsf{P}_t\left(\frac{\pi}{\sin 2\pi t}\left.\frac{(-1)^{n+1}\partial^{n+1}}{\partial t\partial t'\partial t_1\dots\partial t_{n-1}}F_{g-1,n+1}^V(t,t',\mathbf{t})\right|_{t' = t}\right) = \mathsf{P}_t\left(\frac{\pi}{\sin 2\pi t}F_{g-1,n+1}^{EO}(t,-t,\mathbf{t})\right).
    \end{split}
\]
This constitutes the first term on the right-hand side of \eqref{Equ_Proof_1}.

Similarly, the second term on the right-hand side reads as
\[
    \mathsf{P}_t\left(\frac{\pi}{\sin 2\pi t}\sum_{\substack{g_1+g_2 = g,\ \mathbf{t}_1\sqcup \mathbf{t}_2 = \mathbf{t}\\ (g_i,n_i)\neq (0,2)}}F_{g_1,n_1}^{EO}(t,\mathbf{t}_1)F_{g_2,n_2}^{EO}(-t,\mathbf{t}_2)\right).
\]
To continue, note that Eynard and Orantin's notation in \eqref{Equ_Proof_1} allows the occurrence of either $(g_1,n_1) = (0,2)$ or $(g_2,n_2) = (0,2)$, which differs from the notation of our equations. Therefore, the expression above includes every summand in the second term on the right-hand side of \eqref{Equ_Proof_1} except for the two special cases mentioned. These summands correspond to the following partitions of $\mathbf{t}$:
\[
    (\mathbf{t}_1,\mathbf{t}_2) = (\{t_j\},\mathbf{t}\backslash t_j)\text{ or }(\mathbf{t}\backslash t_j,\{t_j\}),\ j=1,\dots, n-1.
\]
They arise from the third term on the right-hand side as follows:
\[
    \begin{split}
        & \sum_{j=1}^{n-1}\frac{(-1)^{n-1}\partial^{n-1}}{\partial t_1\dots\partial t_{n-1}}\mathsf{P}_t\left(\frac{2\pi t_j}{(t^2 - t_j^2)\sin 2\pi t}\frac{\partial}{\partial t}F_{g,n-1}^V(t,\mathbf{t}\backslash t_j)\right) \\
        & = \sum_{j=1}^{n-1}\mathsf{P}_t\left(\frac{\pi}{\sin 2\pi t}\left(\partial_{t_j}\frac{2 t_j}{t^2 - t_j^2}\right)(-1)^{n-1}\left(\partial_{t,\mathbf{t}\backslash t_j}F_{g,n-1}^V(t,\mathbf{t}\backslash t_j)\right)\right) \\
        & = \sum_{j=1}^{n-1}\mathsf{P}_t\left(\frac{\pi}{\sin 2\pi t}\left(F_{0,2}^{EO}(t,t_j)+F_{0,2}^{EO}(-t,t_j)\right)F_{g,n-1}^{EO}(t,\mathbf{t}\backslash t_j)\right) \\
        & = \mathsf{P}_t\left(\frac{\pi}{\sin 2\pi t}\sum_{\substack{(g_1,n_1) = (0,2)\\\text{or}\ (g_2,n_2) = (0,2)}}F_{g_1,n_1}^{EO}(t,\mathbf{t}_1)F_{g_2,n_2}^{EO}(-t,\mathbf{t}_2)\right).
    \end{split}
\]
These equations conclude \eqref{Equ_Proof_1}. Hence, Eynard and Orantin's topological recursion holds.
\end{proof}
One more lemma is required for proving statement (2) in Theorem \ref{thm:1:5}.
\begin{lemma}\label{lem:4:6}
    For any $m\in \mathbb{N}$ and $\omega>0$, the following formula for the Laplace transform holds:
    \begin{equation}\label{equ:3:99}
        \int_0^\infty ((x+\omega\mathrm{i})^m - (x-\omega\mathrm{i})^m)\mathrm{e}^{-tx}dx = \mathsf{P}_t\left(2\mathrm{i}\sin \omega t\int_0^\infty x^m\mathrm{e}^{-tx}dx\right),\ \forall t\in\mathbb{C}_+.
    \end{equation}
\end{lemma}
The proof is included in Appendix \ref{Append:1}.
\begin{proof}[Proof of Theorem \ref{thm:1:5}, statement (2)]
     We will show that Equation \eqref{equ:1:6} can be obtained by multiplying the cut-and join equation \eqref{equ:1:2} with $\mathrm{e}^{-bt -\mathbf{b}\cdot\mathbf{t}}$ and integrating the result by $b$, $b_1$,\dots, $b_{n-1}$ from $0$ to $\infty$ respectively. Utilizing Lemma \ref{lem:4:4} with $\omega = 2\pi$, the left-hand side of the resulting equation becomes
     \[
     \begin{split}
         & \frac{1}{4\pi\mathrm{i}}\int_{\mathbb{R}_+^{n}}((b+2\pi\mathrm{i})V_{g,n}(b+2\pi\mathrm{i},\mathbf{b}) - (b-2\pi\mathrm{i})V_{g,n}(b-2\pi\mathrm{i},\mathbf{b}))\mathrm{e}^{-bt-\mathbf{b}\cdot\mathbf{t}}dbd\mathbf{b} \\
         & = \mathsf{P}_t\left(\frac{\sin 2\pi t}{2\pi}\int_{\mathbb{R}_+^{n}} bV_{g,n}(b,\mathbf{b})\mathrm{e}^{-bt-\mathbf{b}\cdot\mathbf{t}}dbd\mathbf{b}\right) = -\mathsf{P}_t\left(\frac{\sin 2\pi t}{2\pi}F_{g,n}^V(t,\mathbf{t})\right).
     \end{split}
     \]
     Meanwhile, the first term on the right-hand side becomes
     \[
     \begin{split}
         & \frac{1}{2}\int_{\mathbb{R}_+^n}\left(\iint_{\substack{b'+b''\leq b\\ b',b''\geq 0}}V_{g-1,n+1}(b',b'',\mathbf{b})b'b''db'db''\right)\mathrm{e}^{-bt-\mathbf{b}\cdot\mathbf{t}}dbd\mathbf{b} \\
         & = \frac{1}{2}\int_{\mathbb{R}_+^{n+1}}\left(\int_{b'+b''}^\infty \mathrm{e}^{-bt}db\right)V_{g-1,n+1}(b',b'',\mathbf{b})\mathrm{e}^{-\mathbf{b}\cdot\mathbf{t}}b'b''db'db''d\mathbf{b} \\
         & = \frac{1}{2}\int_{\mathbb{R}_+^{n+1}}\frac{\mathrm{e}^{-(b'+b'')t}}{t}V_{g-1,n+1}(b',b'',\mathbf{b})\mathrm{e}^{-\mathbf{b}\cdot\mathbf{t}}b'b''db'db''d\mathbf{b} \\
         & = \frac{1}{2t}\left.\frac{\partial^2}{\partial t\partial t'}F_{g-1,n+1}(t,t',\mathbf{t})\right|_{t' = t}.
     \end{split}
     \]
     Similarly, the second term on the right-hand side becomes
     \[
     \begin{split}
         \frac{1}{2t}\sum_{\substack{g_1+g_2 = g\\ \mathbf{t}_1\sqcup\mathbf{t}_2 = \mathbf{t}}}\left(\frac{\partial}{\partial t}F_{g_1,n_1}^V(t,\mathbf{t}_1)\right)\left(\frac{\partial}{\partial t}F_{g_2,n_2}^V(t,\mathbf{t}_2)\right).
     \end{split}
     \]
     Since $V_{g,n-1}$ is an even function, we replace $(b-b_j)$ with $|b-b_j|$ when computing the third term:
     \[
     \resizebox{\textwidth}{!}{$
     \begin{aligned}
         & \frac{1}{2}\sum_{j=1}^{n-1}\int_{\mathbb{R}_+^n}\left(\int_{0}^{b+b_j}+\int_{0}^{|b-b_j|}V_{g,n-1}(b',\mathbf{b}\backslash b_j) b'db'\right)\mathrm{e}^{-bt-\mathbf{b}\cdot\mathbf{t}}dbd\mathbf{b}\\
         & = \frac{1}{2}\sum_{j=1}^{n-1}\int_{\mathbb{R}_+^{n-2}}\left(\iiint_{\substack{b'\leq b+b_j\\ b,b',b_j\geq 0}}+\iiint_{\substack{b'\leq |b-b_j|\\ b,b',b_j\geq 0}}V_{g,n-1}(b',\mathbf{b}\backslash b_j) b'\mathrm{e}^{-bt-b_jt_j}dbdb'db_j\right)\mathrm{e}^{-\mathbf{b}\cdot\mathbf{t}+b_jt_j}d(\mathbf{b}\backslash b_j)
     \end{aligned}
     $}
     \]
     Note that:
     \[
         b'\leq |b-b_j|\iff b\geq b_j+b'\ \text{or}\ b_j\geq b+b'.
     \]
     Thus, we can swap the order of integrals and simplify the expression as
     \[
     \resizebox{\textwidth}{!}{$
     \begin{aligned}
         & \frac{1}{2}\sum_{j=1}^{n-1}\int_{\mathbb{R}_+^{n-1}}\left(\iint_{\substack{b+b_j\geq b\\ b,b_j\geq 0}}+\iint_{\substack{b\geq b_j+b'\\ b,b_j\geq 0}}+\iint_{\substack{b_j\geq b+b'\\ b,b_j\geq 0}}\mathrm{e}^{-bt-b_jt_j}dbdb_j\right)V_{g,n-1}(b',\mathbf{b}\backslash b_j) b'\mathrm{e}^{-\mathbf{b}\cdot\mathbf{t}+b_jt_j}db'd(\mathbf{b}\backslash b_j)\\
         & = \frac{1}{2}\sum_{j=1}^{n-1}\int_{\mathbb{R}_+^{n-1}}\left(\frac{2t\mathrm{e}^{-b't_j}}{t_j(t^2-t_j^2)}-\frac{2t_j\mathrm{e}^{-b't}}{t(t^2-t_j^2)}\right)V_{g,n-1}(b',\mathbf{b}\backslash b_j) b'\mathrm{e}^{-\mathbf{b}\cdot\mathbf{t}+b_jt_j}db'd(\mathbf{b}\backslash b_j)\\
         & = \sum_{j=1}^{n-1}\frac{1}{t^2 - t_j^2}\left(\frac{t}{t_j}\int_{\mathbb{R}_+^{n-1}}V_{g,n-1}(b',\mathbf{b}\backslash b_j)\mathrm{e}^{-b't_j-\mathbf{b}\cdot\mathbf{t}+b_jt_j}b'db'd(\mathbf{b}\backslash b_j)\right.\\
         & \left.-\frac{t_j}{t}\int_{\mathbb{R}_+^{n-1}}V_{g,n-1}(b',\mathbf{b}\backslash b_j)\mathrm{e}^{-b't-\mathbf{b}\cdot\mathbf{t}+b_jt_j}b'db'd(\mathbf{b}\backslash b_j)\right)\\
         & = \sum_{j=1}^{n-1}\frac{1}{t^2-t_j^2}\left(\frac{t_j}{t}\frac{\partial}{\partial t}F_{g,n-1}^V(t,\mathbf{t}\backslash t_j) - \frac{t}{t_j}\frac{\partial}{\partial t_j}F_{g,n-1}^V(\mathbf{t})\right).
     \end{aligned}
     $}
     \]
     These conclude equation \eqref{equ:1:6}.
\end{proof}
Compared to the interrelation between the cut-and-joint equations \eqref{equ:1:1} and \eqref{equ:1:2}, the topological recursions \eqref{equ:1:5} and \eqref{equ:1:6} exhibit a stronger connection with each other.
\begin{prop}\label{Thm_Mirz_Equiv}
Topological recursion formulas \eqref{equ:1:5} and \eqref{equ:1:6} are equivalent.
\end{prop}
\begin{proof}
Equation \eqref{equ:1:5} implies that the Laurent series of the following function at $t=0$:
\[
\begin{split}
    & t\frac{\partial}{\partial t}F_{g,n}^V(t,\mathbf{t}) + \pi\csc(2\pi t)\left.\frac{\partial^2}{\partial t\partial t'}F_{g-1,n+1}^V(t,t',\mathbf{t})\right|_{t' = t} \\
    & + \sum_{\substack{g_1+g_2 = g\\\mathbf{t}_1\sqcup \mathbf{t}_2 = \mathbf{t}}}\pi\csc(2\pi t)\left(\frac{\partial}{\partial t}F^V_{g_1,n_1}(t,\mathbf{t}_1)\right)\left(\frac{\partial}{\partial t}F^V_{g_2,n_2}(t,\mathbf{t}_2)\right) \\
    & + \sum_{j=1}^{n-1}2\pi t_j(t^2 - t_j^2)^{-1}\csc(2\pi t)\frac{\partial}{\partial t}F_{g,n-1}^V(t,\mathbf{t}\backslash t_j),
\end{split}
\]
has zero principal part. Therefore, for each choice of $\mathbf{t}$, there is a certain neighborhood of the origin in which the function is holomorphic in $t$.

Since $\frac{\sin 2\pi t}{2\pi t}$ is also holomorphic in a neighborhood of $t = 0$, their product:
\[
\begin{split}
    & \frac{\sin 2\pi t}{2\pi}\frac{\partial}{\partial t}F_{g,n}^V(t,\mathbf{t}) + \frac{1}{2t}\left.\frac{\partial^2}{\partial t\partial t'}F_{g-1,n+1}^V(t,t',\mathbf{t})\right|_{t' = t} \\
    & + \sum_{\substack{g_1+g_2 = g\\\mathbf{t}_1\sqcup \mathbf{t}_2 = \mathbf{t}}}\frac{1}{2t}\left(\frac{\partial}{\partial t}F^V_{g_1,n_1}(t,\mathbf{t}_1)\right)\left(\frac{\partial}{\partial t}F^V_{g_2,n_2}(t,\mathbf{t}_2)\right) \\
    & + \sum_{j=1}^{n-1}\frac{t_j}{t(t^2 - t_j^2)}\frac{\partial}{\partial t}F_{g,n-1}^V(t,\mathbf{t}\backslash t_j)
\end{split}
\]
is also holomorphic. Hence, the Laurent expansion of the expression above at $t = 0$ has its principal part vanish. By Lemma \ref{lem:4:4}, the principal part of the last term reads as
\[
    \sum_{j=1}^{n-1}\left(\frac{t_j}{t(t^2 - t_j^2)}\frac{\partial}{\partial t}F_{g,n-1}^V(t,\mathbf{t}\backslash t_j) - \frac{t}{t_j(t^2 - t_j^2)}\frac{\partial}{\partial t_j}F_{g,n-1}^V(t_j,\mathbf{t}\backslash t_j)\right),
\]
which implies equation \eqref{equ:1:6}. The other direction of Proposition \ref{Thm_Mirz_Equiv} is deduced similarly.



\end{proof}
\subsection{The Moduli space of Super Hyperbolic surfaces}
\label{Sec_3_2}
The Laplace transform $F_{g,n}^{\widehat{V}}$ is defined as
\[
F_{g,n}^{\widehat{V}}(\mathbf{t}) = \int_{\mathbb{R}_+^n}\widehat{V}_{g,n}(\mathbf{b})\mathrm{e}^{-\mathbf{b}\cdot\mathbf{t}}d\mathbf{b}.
\]
For example,
\[
    F_{1,1}^{\widehat{V}}(t) = \int_0^\infty\frac{1}{8}\mathrm{e}^{-bt}db = \frac{1}{8 t}. 
\]
Similarly to the case of usual moduli spaces, the following Lemmas are utilized to prove Theorem \ref{thm:1:6}:
\begin{lemma}\label{lem:4:11}
    For $b,t\in \mathbb{C}_+$, the following integrals and series converge and are equal to each other:
    \begin{equation}
        \int_0^\infty \frac{\mathrm{e}^{-tx}dx}{\cosh (b+x)/4} = -2\sum_{k=1}^\infty\frac{(-1)^k\mathrm{e}^{-b(k-1/2)/2}}{t+k/2-1/4}.
    \end{equation}
    \begin{equation}
        \int_0^\infty \frac{\mathrm{e}^{-tx}dx}{\cosh (b-x)/4} = \frac{4\pi \mathrm{e}^{-t b}}{\cos 2\pi t} - 2\sum_{k=1}^\infty\frac{(-1)^k\mathrm{e}^{-b(k-1/2)/2}}{t-k/2+1/4}.
    \end{equation}
\end{lemma}
\begin{lemma}\label{lem:4:12}
    For $b,t_1,t_2\in \mathbb{C}_+$, the following integrals and series converge and are equal to each other:
    \begin{equation}
        \int_0^\infty\int_0^\infty \frac{\mathrm{e}^{-t_1x_1-t_2x_2}dx_1dx_2}{\cosh (b+x_1+x_2)/4} = -2\sum_{k=1}^\infty\frac{(-1)^k\mathrm{e}^{-bk/2}}{(t_1+k/2-1/4)(t_2+k/2-1/4)},
    \end{equation}
    \begin{equation}
    \resizebox{0.9\textwidth}{!}{$
    \begin{aligned}
        \int_0^\infty\int_0^\infty \frac{\mathrm{e}^{-t_1x_1-t_2x_2}dx_1dx_2}{\cosh (b-x_1+x_2)/4} = \frac{4\pi\mathrm{e}^{-t_1b}}{(t_1+t_2)\cos 2\pi t_1}-2\sum_{k=1}^\infty\frac{(-1)^k\mathrm{e}^{-bk/2}}{(t_1-k/2+1/4)(t_2+k/2-1/4)},
    \end{aligned}
    $}
    \end{equation}
    \begin{equation}
    \resizebox{0.9\textwidth}{!}{$
    \begin{aligned}
        \int_0^\infty\int_0^\infty \frac{\mathrm{e}^{-t_1x_1-t_2x_2}dx_1dx_2}{\cosh (b+x_1-x_2)/4} = \frac{4\pi\mathrm{e}^{-t_2b}}{(t_1+t_2)\cos 2\pi t_2}-2\sum_{k=1}^\infty\frac{(-1)^k\mathrm{e}^{-bk/2}}{(t_1+k/2-1/4)(t_2-k/2+1/4)},
    \end{aligned}
    $}
    \end{equation}
    \begin{equation}
    \resizebox{0.9\textwidth}{!}{$
    \begin{aligned}
        \int_0^\infty\int_0^\infty \frac{\mathrm{e}^{-t_1x_1-t_2x_2}dx_1dx_2}{\cosh (b-x_1-x_2)/4} = \frac{4\pi}{t_1-t_2}\left(\frac{\mathrm{e}^{-t_2 b}}{\cos 2\pi t_2}-\frac{\mathrm{e}^{-t_1 b}}{\cos 2\pi t_1}\right)-2\sum_{k=1}^\infty\frac{(-1)^k\mathrm{e}^{-bk/2}}{(t_1-k/2+1/4)(t_2-k/2+1/4)}.
    \end{aligned}
    $}
    \end{equation}
\end{lemma}
\begin{lemma}\label{lem:4:13}
    For any odd number $m\in\mathbb{N}$ and $t$, the following series converge and equal to each other:
    \[
     \sum_{k=1}^\infty\frac{(-1)^k2t(2/(k-1/2))^m}{t^2-(k-1/2)^2/4}= \mathsf{H}_t\left(\frac{2\pi}{t^m\cos 2\pi t}\right).
    \]
\end{lemma}
\begin{proof}[Outline of Proof of Theorem \ref{thm:1:6}, statement (1)]
    Equation \eqref{equ:1:7} can be obtained by multiplying the cut-and-join equation \eqref{equ:1:3} with $\mathrm{e}^{-bt-\mathbf{b}\cdot\mathbf{t}}$, integrating the result over $b$, $b_1$,\dots, $b_{n-1}$ from $0$ to $\infty$ respectively, and canceling the factor $2\pi$. Namely, integrating the left-hand side of \eqref{equ:1:3} yields the left-hand side of \eqref{equ:1:7}. By applying Lemma \ref{lem:4:11} and Lemma \ref{lem:4:13}, integrating the first and the second term on the right-hand side of \eqref{equ:1:3} yields the first and the second term on the right-hand side of \eqref{equ:1:7}, respectively. By applying Lemma \ref{lem:4:12} and Lemma \ref{lem:4:13}, integrating the third term on the right-hand side of \eqref{equ:1:3} yields the third term on the right-hand side of \eqref{equ:1:7}. All these steps are analogous to the corresponding ones in the proof of statement (1) in Theorem \ref{thm:1:5}.
\end{proof}
The following Lemma is necessary for proving the statement (2) of the Theorem:
\begin{lemma}\label{lem:4:14}
    For any odd number $m\in \mathbb{N}$ and $\omega>0$, the following formula of Laplace transform holds:
    \[
    \int_0^\infty ((x+\omega\mathrm{i})^m + (x-\omega\mathrm{i})^m)\mathrm{e}^{-tx}dx = 2\mathsf{P}_t\left(\cos \omega t\int_0^\infty x^m\mathrm{e}^{-tx}dx\right),\ \forall t\in\mathbb{C}_+.
    \]
\end{lemma}
\begin{proof}[Outline of Proof of Theorem \ref{thm:1:6}, statement (2)]
    Equation \eqref{equ:1:8} can be obtained by multiplying the cut-and-join equation \eqref{equ:1:4} with $\mathrm{e}^{-bt-\mathbf{b}\cdot\mathbf{t}}$ and integrating the result by $b$, $b_1$,\dots, $b_{n-1}$ from $0$ to $\infty$ respectively. Utilizing Lemma \ref{lem:4:14}, integrating the left-hand side of \eqref{equ:1:4} yields the left-hand side of \eqref{equ:1:8}. The correspondence of the right-hand sides of equations \eqref{equ:1:4} and \eqref{equ:1:8} is computed straightforwardly, analogously to the proof of statement (2) in Theorem \ref{thm:1:5}.
\end{proof}
Applying Lemma \ref{lem:4:4}, we obtain the following analog of Proposition \ref{Thm_Mirz_Equiv}:
\begin{prop}
Topological recursion formulas \eqref{equ:1:7} and \eqref{equ:1:8} are equivalent.
\end{prop}

%% file: Sec_5.tex
\section{The Top- and Lowest Degree Terms of the Topological Recursions}
\label{Chp_4}
The volume functions $V_{g,n}(\mathbf{b})$ and $\widehat{V}_{g,n}(\mathbf{b})$ are polynomials in the variables $b_1,\dots,b_n$. Due to their polynomial nature, it is convenient to consider the top degree terms of the new recursion formulas \eqref{equ:1:2} and \eqref{equ:1:4}. This section demonstrates that the DVV identity can be recovered by examining these highest degree terms. Moreover, if we take the lowest degree terms of \eqref{equ:1:2} and \eqref{equ:1:4}, we obtain identities of cohomology classes in $H^*(\overline{\mathcal{M}}_{g,n})$.
\subsection{The Top Degree Terms and the DVV Identity}
For any set $\{d_1,\dots,d_n\}$ with $\sum_{i=1}^n d_i = 3g-3+n$, the \emph{top intersection number} of the $\psi$ classes of the compact moduli space $\overline{\mathcal{M}}_{g,n}$ is defined as \cite{mirzakhani2007weil}
\[
\langle\tau_{d_1},\dots,\tau_{d_n}\rangle_g = \int_{\overline{\mathcal{M}}_{g,n}}\prod_{i=1}^n\psi_i^{d_i}.
\]
The top intersection number is connected to the top degree term of the volume $V_{g,n}$ by the formula \eqref{is:1}. The following fact is self-evident from \eqref{is:1}:
\begin{prop}\label{prop:5:1}
    The top degree of the function $V_{g,n}(\mathbf{b})$ is $2(3g-3+n)$. Moreover, for any set $\{d_1,\dots,d_n\}$ with $\sum_{i=1}^n d_i = 3g-3+n$, the coefficient of the $\prod_{i=1}^n b_i^{2d_i}$ term in $V_{g,n}(\mathbf{b})$ is (note that $V_{g,n}$ is even)
    \[
    \frac{\langle\tau_{d_1},\dots,\tau_{d_n}\rangle_g}{\prod_{i=1}^n(2^{d_i})(d_i!)} = \frac{\langle\tau_{d_1},\dots,\tau_{d_n}\rangle_g}{\prod_{i=1}^n(2d_i)!!}.
    \]
\end{prop}
\begin{theorem}\label{thm:5:1}
     The top degree terms of the polynomial cut-and-join equation \eqref{equ:1:2} recover the DVV identity, \cite{dijkgraaf1991loop}:
    \[
    \begin{split}
        & (2d+1)!!\langle\tau_{d},\tau_{\mathbf{d}}\rangle_g = \frac{1}{2}\sum_{\mu+\nu = d-2}(2\mu+1)!!(2\nu+1)!!\langle \tau_\mu,\tau_\nu,\tau_{\mathbf{d}}\rangle_{g-1}\\
        & + \frac{1}{2}\sum_{\substack{g_1+g_2 = g\\ \mathbf{d}_1\sqcup\mathbf{d}_2 = \mathbf{d}}}(2\mu+1)!!(2\nu+1)!!\langle\tau_\mu,\tau_{\mathbf{d}_1}\rangle_{g_1}\langle\tau_\nu,\tau_{\mathbf{d}_2}\rangle_{g_2} + \sum_{j=1}^{n-1}\frac{(2d+2d_j-1)!!}{(2d_j-1)!!}\langle \tau_{d+d_j-1},\tau_{\mathbf{d}\backslash d_j}\rangle_g.
    \end{split}
    \]
    Here $d$ is any natural number, $\mathbf{d} = \{d_1,\dots, d_{n-1}\}$ is any set of natural numbers, with $d+\sum_{i=1}^{n-1} d_i = 3g-3+n$.
\end{theorem}
\begin{proof}
    The highest degree of \eqref{equ:1:2} is $2(3g-3+n)$. For any $d,d_1,\dots,d_{n-1}$ with $d+\sum_{i=1}^{n-1} d_i = 3g-3+n$, we consider the coefficients of the $b^{2d}\prod_{i=1}^{n-1} b_i^{2d_i}$ terms in \eqref{equ:1:2}.

    Note that for any polynomial $P(x)$, the top degree terms of $(P(x+a) - P(x-a))/2a$ and $P'(x)$ are equal. Thus, the left-hand side of \eqref{equ:1:2} has the same top degree term as
    \[
    \frac{\partial}{\partial b}\left(bV_{g,n}(b,\mathbf{b})\right).
    \]
    By Proposition \ref{prop:5:1}, the $b^{2d}\prod_{i=1}^{n-1} b_i^{2d_i}$ term of the left-hand side is
    \[
    \begin{split}
        \frac{\partial}{\partial b}\left(\frac{\langle\tau_d,\tau_{\mathbf{d}}\rangle_g}{(2d)!!\prod_{i=1}^{n-1}(2d_i)!!}b^{2d+1}\prod_{i=1}^{n-1} b_i^{2d_i}\right) = \frac{(2d+1)\langle\tau_d,\tau_{\mathbf{d}}\rangle_g}{(2d)!!\prod_{i=1}^{n-1}(2d_i)!!}b^{2d}\prod_{i=1}^{n-1} b_i^{2d_i}.
    \end{split}
    \]
    Multiple terms in $V_{g-1,n+1}(b',b'',\mathbf{b})$ contribute to the $b^{2d}\prod_{i=1}^{n-1} b_i^{2d_i}$ term on the right-hand side of \eqref{equ:1:2}. Specifically, the $b^{2d}\prod_{i=1}^{n-1} b_i^{2d_i}$ term arises from integrating the ${b'}^{2\mu}{b''}^{2\nu}\prod_{i=1}^{n-1} b_i^{2d_i}$ terms, where $\mu+\nu = d-2$. Integrating these terms yields
    \[
    \begin{split}
        & \frac{1}{2}\iint_{\substack{b',b''\geq 0\\ b'+b''\leq b}}\sum_{\mu+\nu = d-2}\frac{\langle\tau_{\mu},\tau_{\nu},\tau_{\mathbf{d}}\rangle_{g-1}}{(2\mu)!!(2\nu)!!\prod_{i=1}^{n-1}(2d_i)!!}{b'}^{2\mu+1}{b''}^{2\nu+1}\prod_{i=1}^{n-1} b_i^{2d_i}db'db''\\
        & = \frac{1}{2}\sum_{\mu+\nu = d-2}\frac{(2\mu+1)!(2\nu+1)!}{(2\mu+2\nu+4)!}\frac{\langle\tau_{\mu},\tau_{\nu},\tau_{\mathbf{d}}\rangle_{g-1}}{(2\mu)!!(2\nu)!!\prod_{i=1}^{n-1}(2d_i)!!}b^{2\mu+2\nu+4}\prod_{i=1}^{n-1} b_i^{2d_i}\\
        & = \frac{1}{2}\sum_{\mu+\nu = d-2}\frac{(2\mu+1)!!(2\nu+1)!!\langle\tau_{\mu},\tau_{\nu},\tau_{\mathbf{d}}\rangle_{g-1}}{(2d)!\prod_{i=1}^{n-1}(2d_i)!!}b^{2d}\prod_{i=1}^{n-1} b_i^{2d_i}.
    \end{split}
    \]
    The $b^{2d}\prod_{i=1}^{n-1} b_i^{2d_i}$ terms arising from the second part of the right-hand side of \eqref{equ:1:2} are computed similarly:
    \[
        \frac{1}{2}\sum_{\substack{g_1+g_2 = g\\ \mathbf{d}_1\sqcup\mathbf{d}_2 = \mathbf{d}}}\frac{(2\mu+1)!!(2\nu+1)!!\langle\tau_\mu,\tau_{\mathbf{d}_1}\rangle_{g_1}\langle\tau_\nu,\tau_{\mathbf{d}_2}\rangle_{g_2}}{(2d)!\prod_{i=1}^{n-1}(2d_i)!!}b^{2d}\prod_{i=1}^{n-1} b_i^{2d_i}.
    \]
    Here, $d_\mu + \sum_{d_i\in \mathbf{d}_1}d_i = 3g_1-3+n_1$, and $d_\nu + \sum_{d_i\in \mathbf{d}_2}d_i = 3g_2-3+n_2$.

    For the third part of the right-hand side of \eqref{equ:1:2}, we consider the top degree terms in $V_{g,n-1}(b',\mathbf{b}\backslash b_j)$. If for any $i\neq j$, the degree of $b_i$ of such a term is $2d_i$, then the degree of $b'$ must be
    \[
    2(3g-3+(n-1))-\sum_{i\neq j}2d_i = 2(d+d_j-1).
    \]
    Integrating the $b'^{2(d+d_j-1)}\prod_{i\neq j}b_i^{2d_i}$ term in $V_{g,n-1}(b',\mathbf{b}\backslash b_j)$ yields:
    \[
    \begin{split}
        & \frac{1}{2}\left(\int_0^{b+b_j}+\int_0^{b-b_j}\right)\frac{\langle\tau_{d+d_j-1},\tau_{\mathbf{d}\backslash d_j}\rangle_g}{(2(d+d_j-1))!!\prod_{i\neq j}(2d_i)!!}b'^{2(d+d_j-1)+1}\prod_{i\neq j}b_i^{2d_i}db'\\
        & = \frac{1}{2}\frac{(b+b_j)^{2(d+d_j)}+(b-b_j)^{2(d+d_j)}}{2(d+d_j)}\frac{\langle\tau_{d+d_j-1},\tau_{\mathbf{d}\backslash d_j}\rangle_g}{(2(d+d_j-1))!!\prod_{i\neq j}(2d_i)!!}\prod_{i\neq j}b_i^{2d_i}.
    \end{split}
    \]
    The $b^{2d}\prod_{i=1}^{n-1} b_i^{2d_i}$ term in the equation above is
    \[
    \begin{split}
        & \frac{1}{2}\frac{2\binom{2(d+d_j)}{2d}b^{2d}b_j^{2d_j}}{2(d+d_j)}\frac{\langle\tau_{d+d_j-1},\tau_{\mathbf{d}\backslash d_j}\rangle_g}{(2(d+d_j-1))!!\prod_{i\neq j}(2d_i)!!}\prod_{i\neq j}b_i^{2d_i} \\
        & = \frac{(2(d+d_j))!\langle\tau_{d+d_j-1},\tau_{\mathbf{d}\backslash d_j}\rangle_g}{(2d)!(2d_j)!(2(d+d_j))!!\prod_{i\neq j}(2d_i)!!}b^{2d}\prod_{i=1}^{n-1}b_i^{2d_i} \\
        & = \frac{(2(d+d_j)-1)!!\langle\tau_{d+d_j-1},\tau_{\mathbf{d}\backslash d_j}\rangle_g}{(2d)!(2d_j-1)!!\prod_{i=1}^{n-1}(2d_i)!!}b^{2d}\prod_{i=1}^{n-1}b_i^{2d_i}.
    \end{split}
    \]
    We conclude by combining the coefficients of the $b^{2d}\prod_{i=1}^{n-1} b_i^{2d_i}$ terms computed above and multiplying both sides by the factor $(2d)!\prod_{i=1}^{n-1}(2d_i)!!$.
\end{proof}
We turn to the case of super moduli spaces. Similarly, for any set $\{d_1,\dots,d_n\}$ with $\sum_{i=1}^nd_i = g-1$, the top intersection number with the Norbury class is defined as \cite{norbury2020enumerative}
\[
    \langle\tau_{d_1},\dots,\tau_{d_n}\rangle_g^\Theta = \int_{\overline{\mathcal{M}}_{g,n}}\Theta_{g,n}\prod_{i=1}^n\psi_i^{d_i}.
\]
By formula \eqref{is:2}, the following fact is self-evident:
\begin{prop}\label{prop:5:2}
    The top degree of the function $\widehat{V}_{g,n}(\mathbf{b})$ is $2g-2$. Moreover, for any set $\{d_1,\dots,d_n\}$ with $\sum_{i=1}^n d_i = g-1$, the coefficient of the $\prod_{i=1}^n b_i^{2d_i}$ term in $\widehat{V}_{g,n}(\mathbf{b})$ is
    \[
    \frac{\langle\tau_{d_1},\dots,\tau_{d_n}\rangle_g^\Theta}{\prod_{i=1}^n(2^{d_i})(d_i!)} = \frac{\langle\tau_{d_1},\dots,\tau_{d_n}\rangle_g^\Theta}{\prod_{i=1}^n(2d_i)!!}.
    \]
\end{prop}
\begin{theorem}
    The top degree term of the polynomial cut-and-join equation \eqref{equ:1:4} recovers the following:
    \[
    \begin{split}
        & -(2d+1)!!\langle \tau_d,\tau_{\mathbf{d}}\rangle_g^\Theta = \frac{1}{2}\sum_{\mu+\nu = d-1}(2\mu+1)!!(2\nu+1)!!\langle \tau_\mu,\tau_\nu,\tau_{\mathbf{d}}\rangle_{g-1}^\Theta \\
        & + \frac{1}{2}\sum_{\substack{g_1+g_2 = g\\ \mathbf{d}_1\sqcup\mathbf{d}_2 = \mathbf{d}}}(2\mu+1)!!(2\nu+1)!!\langle \tau_\mu,\tau_{\mathbf{d}_1}\rangle_{g_1}^\Theta\langle \tau_\nu,\tau_{\mathbf{d}_2}\rangle_{g_2}^\Theta + \sum_{j=1}^{n-1}\frac{(2d+2d_j+1)!!}{(2d_j-1)!!}\langle \tau_{d+d_j},\tau_{\mathbf{d}\backslash d_j}\rangle_g^\Theta.
    \end{split}
    \]
    Here $d$ is any natural number, $\mathbf{d} = \{d_1,\dots, d_{n-1}\}$ is any set of natural numbers, with $d+\sum_{i=1}^{n-1} d_i = g-1$.
\end{theorem}
\begin{proof}[Outline of the Proof]
    We prove this result by considering the $b^{2d+1}\prod_{i=1}^{n-1}b_i^{2d_i}$ term on both sides of \eqref{equ:1:4} and applying Proposition \ref{prop:5:2}, similarly to the proof of Theorem \ref{thm:5:1}.
\end{proof}
\subsection{The Lowest Degree Terms}
We recover relations for $V_{g,n}$ and $\hat{V}_{g,n}$ by taking the lowest degree terms of the new cut-and-join equations \eqref{equ:1:2} and \eqref{equ:1:4}. These relations involve a complex boundary length $2\pi\mathrm{i}$.

If we set $b=0$ in \eqref{equ:1:2}, the first two terms on the right-hand side vanish. This results in the following relation initially suggested by Do and Norbury:
\begin{prop}[\cite{do2009weil}]
    For any $g,n$ such that $2g+n\geq 3$,
    \[
    V_{g,n}(2\pi\mathrm{i},\mathbf{b}) = \sum_{j=1}^{n-1}\int_0^{b_j}bV_{g,n-1}(b,\mathbf{b}\backslash b_j)db.
    \]
    Here $\mathbf{b} = (b_1,\dots,b_{n-1})$ as before.
\end{prop}
The new cut-and-join equation \eqref{equ:1:2} also implies formulas for higher-order derivatives of $V_{g,n}(2\pi\mathrm{i},\mathbf{b})$:
\begin{theorem}\label{cor:5}
    For any $g,n$ such that $2g+n\geq 3$,
    \begin{equation}\label{equ:5}
    \frac{\partial^2}{\partial b^2}V_{g,n}(2\pi\mathrm{i},\mathbf{b}) = \left(\sum_{j=1}^{n-1}\frac{\partial}{\partial b_j}(b_j V_{g,n-1}(\mathbf{b}))\right) - 2(2g-3+n)V_{g,n-1}(\mathbf{b}).
    \end{equation}
    Here, $\frac{\partial^2}{\partial b^2}$ differentiates with respect to the first variable of $V_{g,n}$.
\end{theorem}
\begin{proof}
We consider the $b^2$ term (including higher-order terms in $b_1,\dots,b_{n-1}$) of both sides of \eqref{equ:1:2}. On the left-hand side, the coefficient is
\[
\begin{split}
    & \frac{1}{4\pi \mathrm{i}}\left.\left(\frac{1}{2!}\frac{\partial^2}{\partial b^2}(bV_{g,n}(b,\mathbf{b}))\right)\right|^{b=2\pi\mathrm{i}}_{b=-2\pi\mathrm{i}}\\
    & = \frac{1}{8\pi \mathrm{i}}\left.\left(b\frac{\partial^2}{\partial b^2}V_{g,n}(b,\mathbf{b})+2\frac{\partial}{\partial b}V_{g,n}(b,\mathbf{b})\right)\right|^{b=2\pi\mathrm{i}}_{b=-2\pi\mathrm{i}}\\
    & = \frac{1}{2}\frac{\partial^2}{\partial b^2}V_{g,n}(2\pi\mathrm{i},\mathbf{b})+ \frac{1}{2\pi\mathrm{i}}\frac{\partial}{\partial b}V_{g,n}(2\pi\mathrm{i},\mathbf{b}).
\end{split}
\]
The first two terms on the right-hand side of \eqref{equ:1:2} do not contribute to the $b^2$ terms. For the third term on the right-hand side, the $b^2$ coefficient reads as
\[
\begin{split}
    & \frac{1}{2}\sum_{j=1}^{n-1}\frac{1}{2}\left(\frac{\partial^2}{\partial b^2}\left.\left(\int_0^b b'V_{g,n-1}(b',\mathbf{b}\backslash b_j)\right)\right|_{b=b_j}+\frac{\partial^2}{\partial b^2}\left.\left(\int_0^b b'V_{g,n-1}(b',\mathbf{b}\backslash b_j)\right)\right|_{b=-b_j}\right)\\
    & = \frac{1}{2}\sum_{j=1}^{n-1}\left.\frac{\partial}{\partial b}(bV_{g,n-1}(b,\mathbf{b}\backslash b_j))\right|_{b=b_j} = \frac{1}{2}\sum_{j=1}^{n-1}\frac{\partial}{\partial b_j}(b_jV_{g,n-1}(\mathbf{b})).
\end{split}
\]
Utilizing Proposition \ref{prop:1:1}:
\[
\frac{1}{2\pi\mathrm{i}}\frac{\partial}{\partial b}V_{g,n}(2\pi\mathrm{i},\mathbf{b}) = (2g-3+n)V_{g,n-1}(\mathbf{b}),
\]
the desired conclusion follows.
\end{proof}
Theorem \ref{cor:5} implies the following relation in terms of the cohomology classes $\kappa_1$, $\psi$ and $\psi_i\in H^*(\mathcal{M}_{g,n})$, $i=1,\dots, n-1$, analogously to \cite{do2009weil}, Lemma 1 and 2.
\begin{Col}
    \[
    \pi_!\left\{\psi^2(\kappa_1-\psi)^{m-1}\prod_{i=1}^{n-1}\psi_i^{d_i}\right\} = \kappa_1^m\prod_{i=1}^{n-1}\psi_i^{d_i},
    \]
    where $\sum_i d_i + m + 1 = 3g-3+n$, $\pi_!: H^*(\overline{\mathcal{M}}_{g,n})\to H^*(\overline{\mathcal{M}}_{g,n-1})$ is the umkehr map induced by the forgetful map $\pi: \overline{\mathcal{M}}_{g,n}\to \overline{\mathcal{M}}_{g,n-1}$.
\end{Col}
\begin{proof}
    The coefficient of $\prod_{i=1}^{n-1}b_i^{2d_i}$ in the left-hand side of \eqref{equ:5} reads as
    \[
    \begin{split}
        & \sum_{j=0}^m\frac{(2j+2)(2j+1)(2\pi\mathrm{i})^{2j}}{2^{\sum_i d_i + j+1}\prod_i d_i! (j+1)!(m-j)!}\int_{\overline{\mathcal{M}}_{g,n}}\psi^{j+1}\prod_i \psi_i^{d_i}(2\pi^2\kappa_1)^{m-j} \\
        & = \frac{2^{m-\sum_i d_i}\pi^{2m}}{\prod_i d_i! m!}\sum_{j=0}^m(2j+1)(-1)^j\binom{m}{j}\int_{\overline{\mathcal{M}}_{g,n}}\psi^{j+1}\prod_i \psi_i^{d_i}\kappa_1^{m-j} \\
        & = \frac{2^{m-\sum_i d_i}\pi^{2m}}{\prod_i d_i! m!}\int_{\overline{\mathcal{M}}_{g,n}}\psi(\kappa_1-\psi)^{m-1}(\kappa_1-(2m+1)\psi)\prod_i\psi_i^{d_i}.
    \end{split}
    \]
    The coefficient of $\prod_{i=1}^{n-1}b_i^{2d_i}$ in the first term on the right-hand side of \eqref{equ:5} is
    \[
    \begin{split}
        & \sum_{k=1}^{n-1}\frac{2d_k+1}{2^{\sum_i d_i}\prod_id_i!m!}\int_{\overline{\mathcal{M}}_{g,n-1}}\prod_{i=1}^{n-1}\psi_i^{d_i}(2\pi^2\kappa_1)^m \\
        & = \sum_{k=1}^{n-1}\frac{2^{m-\sum_i d_i}\pi^{2m}(2d_k+1)}{\prod_id_i!m!}\int_{\overline{\mathcal{M}}_{g,n-1}}\prod_{i=1}^{n-1}\psi_i^{d_i}\kappa_1^m \\
        & = (3(2g-3+n)-2m)\frac{2^{m-\sum_i d_i}\pi^{2m}}{\prod_id_i!m!}\int_{\overline{\mathcal{M}}_{g,n-1}}\prod_{i=1}^{n-1}\psi_i^{d_i}\kappa_1^m.
    \end{split}
    \]
    In the second term, the coefficient is, \cite{do2009weil}
    \[
        -2(2g-3+n)\frac{2^{m-\sum_i d_i}\pi^{2m}}{\prod_id_i!m!}\int_{\overline{\mathcal{M}}_{g,n-1}}\prod_{i=1}^{n-1}\psi_i^{d_i}\kappa_1^m.
    \]
    Thus,
    \[
    \int_{\overline{\mathcal{M}}_{g,n}}\psi(\kappa_1-\psi)^{m-1}(\kappa_1-(2m+1)\psi)\prod_{i=1}^{n-1}\psi_i^{d_i} = ((2g-3+n)-2m)\int_{\overline{\mathcal{M}}_{g,n-1}}\prod_{i=1}^{n-1}\psi_i^{d_i}\kappa_1^m.
    \]
    Subtracting this from the equation shown in \cite{do2009weil}
    \[
    \int_{\overline{\mathcal{M}}_{g,n}}\psi(\kappa_1-\psi)^{m-1}(\kappa_1-\psi)\prod_{i=1}^{n-1}\psi_i^{d_i} = (2g-3+n)\int_{\overline{\mathcal{M}}_{g,n-1}}\prod_{i=1}^{n-1}\psi_i^{d_i}\kappa_1^m,
    \]
    we derive
    \[
    \int_{\overline{\mathcal{M}}_{g,n}}\psi^2(\kappa_1-\psi)^{m-1}\prod_{i=1}^{n-1}\psi_i^{d_i} = \int_{\overline{\mathcal{M}}_{g,n-1}}\prod_{i=1}^{n-1}\psi_i^{d_i}\kappa_1^m.
    \]
\end{proof}
Utilizing Proposition \ref{prop:1:2}, we can similarly show that the cut-and-join equation \eqref{equ:1:4} implies:
\begin{theorem}
    For any $g,n\geq 1$,
    \begin{equation}\label{equ:5:2}
        -2\pi\mathrm{i}\frac{\partial}{\partial b}\widehat{V}_{g,n}(2\pi\mathrm{i},\mathbf{b}) = -\left(\sum_{j=1}^{n-1}\frac{\partial}{\partial b_j}(b_j \widehat{V}_{g,n-1}(\mathbf{b}))\right) +(2g-3+n)\widehat{V}_{g,n-1}(\mathbf{b}).
    \end{equation}
\end{theorem}
This implies the following relation in terms of the cohomology classes $\kappa_1$, $\psi$, $\psi_i$, $i=1,\dots, n-1$, and $\Theta_{g,n}$:
\begin{Col}
    \[
    \pi_!\left\{\psi(\kappa_1-\psi)^m\prod_{i=1}^{n-1}\psi_i^{d_i}\Theta_{g,n}\right\} = \kappa_1^{m+1}\prod_{i=1}^{n-1}\psi_i^{d_i}\Theta_{g,n-1},
    \]
    where $\sum_i d_i + m + 1 = g-1$.
\end{Col}
\begin{proof}
    The coefficient of $\prod_{i=1}^{n-1}b_i^{2d_i}$ in the left-hand side of \eqref{equ:5:2} is
    \[
    \begin{split}
        & \sum_{j=0}^m\frac{-(2j+2)(2\pi\mathrm{i})^{2j+2}}{2^{\sum_i d_i + j + 1}\prod_i d_i! (j+1)!(m-j)!}\int_{\overline{\mathcal{M}}_{g,n}}\psi^{j+1}\prod_{i}\psi_i^{d_i}(2\pi^2\kappa_1)^{m-j}\Theta_{g,n}\\
        & = \sum_{j=0}^m\frac{2^{m+2-\sum_i d_i}\pi^{2m+2}}{\prod_i d_i!m!}\int_{\overline{\mathcal{M}}_{g,n}}(-1)^j\binom{m}{j}\psi^{j+1}\prod_{i}\psi_i^{d_i}\kappa_1^{m-j}\Theta_{g,n} \\
        & = \frac{2^{m+2-\sum_i d_i}\pi^{2m+2}}{\prod_i d_i!m!}\int_{\overline{\mathcal{M}}_{g,n}}\psi(\kappa_1 - \psi)^m\prod_{i}\psi_i^{d_i}\Theta_{g,n}.
    \end{split}
    \]
    The coefficient of $\prod_{i=1}^{n-1}b_i^{2d_i}$ in the first term of the right-hand side is
    \[
    \begin{split}
        & -\sum_{k=1}^{n-1}\frac{2d_k+1}{2^{\sum_i d_i}\prod_id_i!(m+1)!}\int_{\overline{\mathcal{M}}_{g,n-1}}\prod_{i=1}^{n-1}\psi_i^{d_i}(2\pi^2\kappa_1)^{m+1}\Theta_{g,n-1} \\
        & = -\sum_{k=1}^{n-1}\frac{2^{m+1-\sum_i d_i}\pi^{2m+2}(2d_k+1)}{\prod_id_i!(m+1)!}\int_{\overline{\mathcal{M}}_{g,n-1}}\prod_{i=1}^{n-1}\psi_i^{d_i}\kappa_1^{m+1}\Theta_{g,n-1} \\
        & = -(2g-5+n-2m)\frac{2^{m+1-\sum_i d_i}\pi^{2m+2}}{\prod_id_i!(m+1)!}\int_{\overline{\mathcal{M}}_{g,n-1}}\prod_{i=1}^{n-1}\psi_i^{d_i}\kappa_1^{m+1}\Theta_{g,n-1}.
    \end{split}
    \]
    The coefficient of $\prod_{i=1}^{n-1}b_i^{2d_i}$ in the second term is
    \[
    (2g-3+n)\frac{2^{m+1-\sum_i d_i}\pi^{2m+2}}{\prod_id_i!(m+1)!}\int_{\overline{\mathcal{M}}_{g,n-1}}\prod_{i=1}^{n-1}\psi_i^{d_i}\kappa_1^{m+1}\Theta_{g,n-1}.
    \]
    Summing the two terms gives
    \[
    \begin{split}
        & (2m+2)\frac{2^{m+1-\sum_i d_i}\pi^{2m+2}}{\prod_id_i!(m+1)!}\int_{\overline{\mathcal{M}}_{g,n-1}}\prod_{i=1}^{n-1}\psi_i^{d_i}\kappa_1^{m+1}\Theta_{g,n-1} \\
        & = \frac{2^{m+2-\sum_i d_i}\pi^{2m+2}}{\prod_id_i!m!}\int_{\overline{\mathcal{M}}_{g,n-1}}\prod_{i=1}^{n-1}\psi_i^{d_i}\kappa_1^{m+1}\Theta_{g,n-1}.
    \end{split}
    \]
    Therefore,
    \[
    \int_{\overline{\mathcal{M}}_{g,n}}\psi(\kappa_1 - \psi)^m\prod_{i}\psi_i^{d_i}\Theta_{g,n} = \int_{\overline{\mathcal{M}}_{g,n-1}}\prod_{i=1}^{n-1}\psi_i^{d_i}\kappa_1^{m+1}\Theta_{g,n-1}.
    \]
\end{proof}

%% file: Sec_6.tex
\section{The Partition Function of the Highest-degree Terms}
Consider the lowest degree term of $F_{g,n}^V(\mathbf{t})$, denote it by $G_{g,n}^V(\mathbf{t})$. The term is the Laplace transform of the highest degree term of $V_{g,n}(\mathbf{b})$. Similarly, $G_{g,n}^{\widehat{V}}(\mathbf{t})$, the lowest degree term of $F_{g,n}^{\widehat{V}}(\mathbf{t})$, is the Laplace transform of the highest degree term of $\widehat{V}_{g,n}(\mathbf{b})$. The degree $-3(2g-2+n)-1$ part of \eqref{equ:1:6} eliminates the $\mathsf{P}_t$ notation, and shows that
\begin{equation}\label{equ:6:1}
\begin{split}
    & -\frac{\partial}{\partial t}G^V_{g,n}(t,\mathbf{t}) = \frac{1}{2t^2}\left.\frac{\partial^2}{\partial t\partial t'}G^V_{g-1,n+1}(t,t',\mathbf{t})\right|_{t'=t}\\
    & + \frac{1}{2t^2}\sum_{\substack{g_1+g_2 = g\\ \mathbf{t}_1\sqcup\mathbf{t}_2 = \mathbf{t}}}\left(\frac{\partial}{\partial t}G^V_{g_1,n_1}(t,\mathbf{t}_1)\right)\left(\frac{\partial}{\partial t}G^V_{g_2,n_2}(t,\mathbf{t}_2)\right)\\
    & + \sum_{j=1}^{n-1}\frac{t_j}{t^2-t_j^2}\left(\frac{1}{t^2}\frac{\partial}{\partial t}G^V_{g,n-1}(t,\mathbf{t}\backslash t_j) - \frac{1}{t_j^2}\frac{\partial}{\partial t_j}G^V_{g,n-1}(\mathbf{t})\right).
\end{split}
\end{equation}
Similarly, the degree $-(2g-2+n)-1$ part of \eqref{equ:1:8} shows that
\begin{equation}\label{equ:6:2}
\begin{split}
    & -\frac{\partial}{\partial t}G^{\widehat{V}}_{g,n}(t,\mathbf{t}) = \frac{1}{2}\left.\frac{\partial^2}{\partial t\partial t'}G^{\widehat{V}}_{g-1,n+1}(t,t',\mathbf{t})\right|_{t'=t}\\
    & + \frac{1}{2}\sum_{\substack{g_1+g_2 = g\\ \mathbf{t}_1\sqcup\mathbf{t}_2 = \mathbf{t}}}\left(\frac{\partial}{\partial t}G^{\widehat{V}}_{g_1,n_1}(t,\mathbf{t}_1)\right)\left(\frac{\partial}{\partial t}G^{\widehat{V}}_{g_2,n_2}(t,\mathbf{t}_2)\right)\\
    & + \sum_{j=1}^{n-1}\frac{t_j}{t^2-t_j^2}\left(\frac{\partial}{\partial t}G^{\widehat{V}}_{g,n-1}(t,\mathbf{t}\backslash t_j) - \frac{\partial}{\partial t_j}G^{\widehat{V}}_{g,n-1}(\mathbf{t})\right).
\end{split}
\end{equation}
Define the \textbf{partition function} for the volumes $V_{g,n}$ and $\widehat{V}_{g,n}$ by
\[
Z^V(t,\hbar) = \exp\left(\sum_{g=0}^\infty\sum_{n=1}^\infty\frac{1}{n!}\hbar^{2g-2+n}G^V_{g,n}(t,\dots,t)\right),
\]
and
\[
Z^{\widehat{V}}(t,\hbar) = \exp\left(\sum_{g=1}^\infty\sum_{n=1}^\infty\frac{1}{n!}\hbar^{2g-2+n}G^{\widehat{V}}_{g,n}(t,\dots,t)\right).
\]
Here, we additionally define that
\[
G^V_{0,2}(t_1,t_2) = -\log(t_1+t_2),
\]
while $G_{0,1}^V(t)$ is regarded as zero. The main result of this section is that these partition functions satisfy the following differential equations:
\begin{theorem}\label{thm:6:1}
    Set $x = \frac{t^2}{2}$, then the partition function $Z^V(x,\hbar)$ satisfies a differential equation:
    \[
    \left(\frac{\hbar}{2}\frac{\partial^2}{\partial x^2}+\sqrt{2x}\frac{\partial}{\partial x}+\frac{1}{2\sqrt{2x}}\right)Z^V(x,\hbar) = 0.
    \]
\end{theorem}
\begin{theorem}\label{thm:6:2}
    The partition function $Z^{\widehat{V}}(t,\hbar)$ satisfies a differential equation:
    \[
    \left(\frac{\hbar}{2}\frac{\partial^2}{\partial t^2}+ \frac{\partial}{\partial t}+\frac{\hbar}{8t^2}\right)Z^{\widehat{V}}(t,\hbar) = 0.
    \]
\end{theorem}
\subsection{Proof of Theorem \ref{thm:6:1}}
We begin by taking the summation of functions $G^V_{g,n}(t,\dots,t)$ for $(2g-2+n)$ fixed:
\begin{lemma}
    Denote that
    \[
    S_m(t) = \sum_{2g-2+n = m-1}\frac{G^V_{g,n}(t,\dots,t)}{n!}.
    \]
    Then for any $m\geq 2$,
    \begin{equation}\label{equ:6:3}
    \frac{d}{dt}S_{m+1}(t) = -\frac{1}{2t^2}\left(\frac{d^2}{dt^2}S_m(t) +\sum_{\substack{a,b\geq 2\\a+b=m+1}}\frac{dS_a}{dt}(t)\frac{dS_b}{dt}(t)\right)+\frac{1}{t^3}\frac{d}{dt}S_m(t).
    \end{equation}
\end{lemma}
\begin{proof}
    Equation \eqref{equ:6:3} is derived by multiplying \eqref{equ:6:1} with $\frac{1}{(n-1)!}$ and taking the sum over all pairs $(g,n)$ such that $2g-2+n = m$. The proof is similar to Appendix A in \cite{mulase2012spectral}. Specifically, the first and third items on the right-hand side of \eqref{equ:6:1} together provide
    \[
    \begin{split}
        & -\frac{1}{2}\frac{d}{dt}\left(\frac{1}{t^2}\frac{d}{dt}S_m(t)\right) = -\frac{1}{2}\left(-\frac{2}{t^3}\frac{dS_m(t)}{dt}+ \frac{1}{t^2}\frac{d^2S_m(t)}{dt^2}\right)\\
        & = -\frac{1}{2t^2}\frac{d^2S_m(t)}{dt^2} + \frac{1}{t^3}\frac{dS_m(t)}{dt},
    \end{split}
    \]
    and the second item provides
    \[
    -\frac{1}{2t^2}\left(\sum_{\substack{a,b\geq 2\\a+b=m+1}}\frac{dS_a}{dt}(t)\frac{dS_b}{dt}(t)\right).
    \]
\end{proof}
By setting $x = \frac{t^2}{2}$, equation \eqref{equ:6:3} implies that
\begin{equation}\label{equ:6:4}
\sqrt{2x}\frac{d}{dx}S_{m+1}(x) = -\frac{1}{2}\left(\frac{d^2}{dx^2}S_m(x) +\sum_{\substack{a,b\geq 2\\a+b=m+1}}\frac{dS_a}{dx}(x)\frac{dS_b}{dx}(x)\right)+\frac{1}{4x}\frac{d}{dx}S_m(x).
\end{equation}
\begin{lemma}\label{lem:6:2}
    Denote that
    \[
    F(x,\hbar) = \sum_{m=1}^\infty \hbar^{m-1}S_m(x).
    \]
    Then $F(x,\hbar)$ satisfies the differential equation:
    \[
    \sqrt{2x}\frac{dF}{dx}+\frac{\hbar}{2}\frac{d^2F}{dx^2}+\frac{\hbar}{2}\left(\frac{dF}{dx}\right)^2 + \frac{1}{2\sqrt{2x}} = 0.
    \]
\end{lemma}
\begin{proof}
    Note that
    \[
    S_1(t) = \frac{1}{2}G_{0,2}^V(t,t) = -\frac{1}{2}\log(2t) = -\frac{1}{4}\log(4t^2),
    \]
    thus,
    \[
    S_1(x) = -\frac{1}{4}\log(8x),
    \]
    and
    \[
    \frac{d}{dx}S_1(x) = -\frac{1}{4x}.
    \]
    Moreover,
    \[
    S_2(t) = \frac{G_{0,3}^V(t,t,t)}{6} + G_{1,1}^V(t) = \frac{5}{24t^3}.
    \]
    Thus,
    \[
    S_2(x) = \frac{5\sqrt{2}}{96}x^{-3/2},\ \frac{d}{dx}S_2(x) = -\frac{5\sqrt{2}}{64}x^{-5/2}.
    \]
    By multiplying equation \eqref{equ:6:4} with $\hbar^m$ and taking the summation for $m\geq 2$, we obtain
    \[
    \begin{split}
        & \sqrt{2x}\frac{d}{dx}\left(\sum_{m=3}^\infty \hbar^{m-1}S_m(x)\right) = -\frac{\hbar}{2}\frac{d^2}{dx^2}\left(\sum_{m=2}^\infty \hbar^{m-1}S_m(x)\right) \\
        & -\frac{\hbar}{2}\sum_{m\geq 3}\left(\sum_{\substack{a,b\geq 2\\ a+b=m+1}}\frac{\hbar^{a-1} d S_a(x)}{d x}\frac{\hbar^{b-1} d S_b(x)}{d x}\right) - \hbar \sum_{m=2}^\infty \frac{dS_1(x)}{dx}\frac{d S_m(x)}{dx}.
    \end{split}
    \]
    That is,
    \[
    \begin{split}
        & \sqrt{2x}\frac{d}{dx}(F(x,\hbar) - \hbar S_2(x) - S_1(x)) + \frac{\hbar}{2}\frac{d^2}{dx^2}\left(F(x,\hbar) - S_1(x)\right)\\
        & + \frac{\hbar}{2}\left(\left(\frac{d}{dx}F(x,\hbar)\right)^2 - \left(\frac{d}{dx}S_1(x)\right)^2\right) = 0,
    \end{split}
    \]
    or
    \[
        \sqrt{2x}\frac{dF}{dx}+\frac{\hbar}{2}\frac{d^2F}{dx^2}+\frac{\hbar}{2}\left(\frac{dF}{dx}\right)^2 + \frac{1}{2\sqrt{2x}} = \sqrt{2x}\hbar\frac{dS_2}{dx} + \frac{\hbar}{2}\frac{d^2 S_1}{dx^2}+\frac{\hbar}{2}\left(\frac{dS_1}{dx}\right)^2.
    \]
    The right-hand side is computed as
    \[-\sqrt{2}x^{1/2}\hbar\frac{(5\sqrt{2})}{64}x^{-5/2} + \frac{\hbar}{2}\frac{1}{4x^2} + \frac{\hbar}{2}\left(-\frac{1}{4x}\right)^2
        = -\frac{5\hbar}{32}x^{-2} + \frac{\hbar}{8}x^{-2}+\frac{\hbar}{32}x^{-2} = 0.
    \]
    The conclusion follows.
\end{proof}
\begin{proof}[Proof of Theorem \ref{thm:6:1}]
    Note that $Z^V(x,\hbar) = \exp F(x,\hbar)$. The conclusion of Theorem \ref{thm:6:1} is straightforward from Lemma \ref{lem:6:2} via a variable change.
\end{proof}
\subsection{Proof of Theorem \ref{thm:6:2}}
\begin{lemma}\label{lem:6:4}
    Denote that
    \[
    S_m(t) = \sum_{2g-2+n = m-1}\frac{G^{\widehat{V}}_{g,n}(t,\dots,t)}{n!}.
    \]
    Then for any $m\geq 2$,
    \begin{equation}\label{equ:6:5}
        \frac{d}{dt}S_{m+1}(t) = -\frac{1}{2}\left(\frac{d^2}{dt^2}S_m(t) +\sum_{\substack{a,b\geq 2\\a+b=m+1}}\frac{dS_a}{dt}(t)\frac{dS_b}{dt}(t)\right).
    \end{equation}
\end{lemma}
\begin{lemma}
    Denote that
    \[
    F(t,\hbar) = \sum_{m=2}^\infty \hbar^{m-1}S_m(t).
    \]
    Then $F(x,\hbar)$ satisfies the differential equation:
    \[
    \frac{dF}{dt}+\frac{\hbar}{2}\frac{d^2F}{dt^2} + \frac{\hbar}{2}\left(\frac{dF}{dt}\right)^2 + \frac{\hbar}{8t^2} = 0.
    \]
\end{lemma}
\begin{proof}
    By multiplying \eqref{equ:6:5} with $\hbar^m$ and taking the summation for $m\geq 2$, we obtain,
    \[
    \frac{d}{dt}\left(F(t,\hbar) - \hbar S_2(t)\right)+ \frac{\hbar}{2}\frac{d^2F}{dt^2}(t,\hbar) + \frac{\hbar}{2}\left(\frac{dF}{dt}(t,\hbar)\right)^2 = 0.
    \]
    Since
    \[
    \hbar\frac{dS_2}{dt} = \hbar\frac{d}{dt}\left(\frac{1}{8t}\right) = -\frac{\hbar}{8t^2},
    \]
    the conclusion follows.
\end{proof}
\begin{proof}[Proof of Theorem \ref{thm:6:2}]
    Note that $Z^{\widehat{V}}(t,\hbar) = \exp F(t,\hbar)$. The conclusion of Theorem \ref{thm:6:2} is straightforward from Lemma \ref{lem:6:4} via a variable change.
\end{proof}

%% file: Sec_Append.tex
\section{Appendix: Proofs of Lemmas in Section \ref{Chp_3}}\label{Append:1}
\renewcommand{\thesubsection}{\Alph{subsection}}

\begin{proof}[Proof of Lemma \ref{lem:4:1}]
We start by proving the formula \eqref{equ:3:1}. We rewrite the integrand of the left-hand side using a series expansion:
\[
\frac{\mathrm{e}^{-tx}}{1+\mathrm{e}^{(b+x)/2}} = \mathrm{e}^{-tx}\left(\sum_{k=1}^{\infty}(-1)^{k-1}\mathrm{e}^{-(b+x)k/2}\right).
\]
Since $b,t\in\mathbb{C}_+$, the series converges absolutely, allowing us to interchange the sum and the integral. Thus,
\[
\begin{split}
    & \int_0^\infty \frac{\mathrm{e}^{-tx}dx}{1+\mathrm{e}^{(b+x)/2}} = \int_0^\infty \mathrm{e}^{-tx}\left(\sum_{k=1}^{\infty}(-1)^{k-1}\mathrm{e}^{-(b+x)k/2}\right)dx \\
    & = \sum_{k=1}^\infty (-1)^{k-1}\left(\int_0^\infty \mathrm{e}^{-tx}\mathrm{e}^{-(b+x)k/2} dx\right) = -\sum_{k=1}^\infty \frac{(-1)^k\mathrm{e}^{-bk/2}}{t+k/2}.
\end{split}
\]

Next, we prove that the formula \eqref{equ:3:2} holds for $\frac{n}{2}<\mathrm{Re}(t)<\frac{n+1}{2}$ using induction on $n$. 
\textbf{Base case} ($n=0$):
We start by splitting the integral into two parts:
\[
    \int_0^\infty \frac{\mathrm{e}^{-tx}dx}{1+\mathrm{e}^{(b-x)/2}} = \int_{-\infty}^\infty \frac{\mathrm{e}^{-tx}dx}{1+\mathrm{e}^{(b-x)/2}} - \int_{-\infty}^0 \frac{\mathrm{e}^{-tx}dx}{1+\mathrm{e}^{(b-x)/2}}.
\]
Note that the condition $0<\mathrm{Re}(t)<\frac{1}{2}$ guarantees the convergence of the integral for both $x\to \infty$ and $x\to -\infty$. For the first term:
\[
\begin{split}
    & \int_{-\infty}^\infty \frac{\mathrm{e}^{-tx}dx}{1+\mathrm{e}^{(b-x)/2}} \overset{u = (b-x)/2}{=\joinrel=\joinrel=\joinrel=} 2\mathrm{e}^{-tb}\int_{-\infty}^\infty\frac{\mathrm{e}^{2tu}}{1+\mathrm{e}^u}du \\
    & = \frac{2\mathrm{e}^{-tb}}{1-\mathrm{e}^{4\pi\mathrm{i}t}}\left(\int_{-\infty}^\infty - \int_{2\pi \mathrm{i}-\infty}^{2\pi \mathrm{i} + \infty}\right)\frac{\mathrm{e}^{2tu}}{1+\mathrm{e}^u}du = \frac{2\mathrm{e}^{-tb}}{1-\mathrm{e}^{4\pi\mathrm{i}t}}\lim_{r\to\infty}\int_{C_r}\frac{\mathrm{e}^{2tu}}{1+\mathrm{e}^u}du,
\end{split}
\]
where $C_r$ denotes the rectangular contour given by
\[
    -r\to r\to r+2\pi \mathrm{i}\to-r + 2\pi\mathrm{i}\to -r,
\]
shown in Figure \ref{fig:2}.
\begin{figure}[H]
    \centering
    \begin{tikzpicture}[scale=0.20]
        \draw[thick,->] (-20,0) -- (0,0);
        \draw[thick,-] (0,0) -- (20,0) node[anchor=west]{$r$};
        \draw[thick,->] (20,0) -- (20,3.14);
        \draw[thick,-] (20,3.14) -- (20,6.28) node[anchor=west]{$r+2\pi\mathrm{i}$};
        \draw[thick,->] (20,6.28) -- (0,6.28);
        \draw[thick,-] (0,6.28) -- (-20,6.28) node[anchor=east]{$-r+2\pi\mathrm{i}$};
        \draw[thick,->] (-20,6.28) -- (-20,3.14);
        \draw[thick,-] (-20,3.14) -- (-20,0) node[anchor=east]{$-r$};
    \end{tikzpicture}
    \caption{The path $C_r$.}
    \label{fig:2}
\end{figure}
Since the integrand has a singularity at $u = \pi\mathrm{i}$ inside $C_r$, the residue theorem implies:
\[
    \frac{2\mathrm{e}^{-tb}}{1-\mathrm{e}^{4\pi\mathrm{i}t}}\lim_{r\to\infty}\int_{C_r}\frac{\mathrm{e}^{2tu}}{1+\mathrm{e}^u}du = \frac{(-2\pi\mathrm{i})2\mathrm{e}^{-tb}}{1-\mathrm{e}^{4\pi\mathrm{i}t}}\mathrm{Res}_{u = \pi\mathrm{i}} \frac{\mathrm{e}^{2tu}}{1+\mathrm{e}^u} = \frac{2\pi \mathrm{e}^{-tb}}{\sin 2\pi t}.
\]
For the second term:
\[
    - \int_{-\infty}^0 \frac{\mathrm{e}^{-tx}dx}{1+\mathrm{e}^{(b-x)/2}} = -\sum_{k=1}^\infty \frac{(-1)^k\mathrm{e}^{-bk/2}}{t-k/2},
\]
similarly to \eqref{equ:3:1}. Combining these results concludes the case $n=0$. 

\textbf{Induction}: Assuming that \eqref{equ:3:2} holds for $\frac{n-1}{2}<\mathrm{Re}(t)<\frac{n}{2}$. We will prove it for $\frac{n}{2}<\mathrm{Re}(t)<\frac{n+1}{2}$. Indeed,
\[
\begin{split}
    & \int_0^\infty \frac{\mathrm{e}^{-tx}dx}{1+\mathrm{e}^{(b-x)/2}} = \mathrm{e}^{-b/2}\int_0^\infty \frac{\mathrm{e}^{-(t-1/2)x}\mathrm{e}^{(b-x)/2}dx}{1+\mathrm{e}^{(b-x)/2}}\\
    & = \mathrm{e}^{-b/2}\left(\int_0^\infty \mathrm{e}^{-(t-1/2)x} dx - \int_0^\infty \frac{\mathrm{e}^{-(t-1/2)x}dx}{1+\mathrm{e}^{(b-x)/2}}\right).
\end{split}
\]
Applying the induction hypothesis to \eqref{equ:3:2} with $t$ being replaced by $t - 1/2$:
\[
\begin{split}
    & \mathrm{e}^{-b/2}\left(\int_0^\infty \mathrm{e}^{-(t-1/2)x} dx - \int_0^\infty \frac{\mathrm{e}^{-(t-1/2)x}dx}{1+\mathrm{e}^{(b-x)/2}}\right) \\
    & = \mathrm{e}^{-b/2}\left(\frac{1}{t-1/2} - \frac{2\pi\mathrm{e}^{-(t-1/2)b}}{\sin 2\pi (t-1/2)} + \sum_{k=1}^\infty \frac{(-1)^k\mathrm{e}^{-bk/2}}{t-1/2 -k/2}\right) \\
    & = \frac{2\pi \mathrm{e}^{-tb}}{\sin 2\pi t} + \frac{\mathrm{e}^{-b/2}}{t-1/2} +\sum_{k=1}^\infty \frac{(-1)^k\mathrm{e}^{-b(k+1)/2}}{t-(k+1)/2} = \frac{2\pi \mathrm{e}^{-tb}}{\sin 2\pi t}-\sum_{k=1}^\infty \frac{(-1)^k\mathrm{e}^{-bk/2}}{t-k/2}.
\end{split}
\]
By induction, the formula \eqref{equ:3:2} holds for $b,t\in\mathbb{C}_+$ with $\mathrm{Re}(t) \notin \frac{1}{2}\mathbb{Z}$. 

Finally, the result for $\mathrm{Re}(t) \in \frac{1}{2}\mathbb{Z}$ follows by continuity.

\end{proof}
\begin{proof}[Proof of Lemma \ref{lem:4:2}]
    The proof of equation \eqref{equ:3:3} is straightforward and similar to the proof of \eqref{equ:3:1}, while the proofs of equations \eqref{equ:3:4} and \eqref{equ:3:5} follow the same pattern as \eqref{equ:3:2}. Hence, we focus on proving equation \eqref{equ:3:6}. 
    
    For any $\frac{n_1}{2}<\mathrm{Re}(t_1)<\frac{n_1+1}{2}$ and $\frac{n_2}{2}<\mathrm{Re}(t_2)<\frac{n_2+1}{2}$, we prove formula \eqref{equ:3:6} by double induction on $n_1$ and $n_2$. For the base case where $(n_1,n_2) = (0,0)$, we assume that $\mathrm{Re}(t_2)>\mathrm{Re}(t_1)$. Thus,
    \[
        \int_0^\infty\int_0^\infty\frac{\mathrm{e}^{-t_1x_1-t_2x_2}dx_1dx_2}{1+\mathrm{e}^{(b-x_1-x_2)/2}} = \int_0^\infty\int_{-\infty}^\infty\frac{\mathrm{e}^{-t_1x_1-t_2x_2}dx_1dx_2}{1+\mathrm{e}^{(b-x_1-x_2)/2}} - \int_0^\infty\int_{-\infty}^0\frac{\mathrm{e}^{-t_1x_1-t_2x_2}dx_1dx_2}{1+\mathrm{e}^{(b-x_1-x_2)/2}}.
    \]
    The first term is computed similarly to the proof of \eqref{equ:3:2}:
    \[
        \int_0^\infty\int_{-\infty}^\infty\frac{\mathrm{e}^{-t_1x_1-t_2x_2}dx_1dx_2}{1+\mathrm{e}^{(b-x_1-x_2)/2}} = \int_0^\infty \mathrm{e}^{-t_2x_2}\frac{2\pi\mathrm{e}^{-t_1(b-x_2)}dx_2}{\sin 2\pi t_1} = -\frac{1}{t_1 - t_2}\frac{2\pi\mathrm{e}^{-t_1 b}}{\sin 2\pi t_1},
    \]
    where the convergence follows from the assumption $\mathrm{Re}(t_2)>\mathrm{Re}(t_1)$. The second term is computed similarly to \eqref{equ:3:5}:
    \[
        - \int_0^\infty\int_{-\infty}^0\frac{\mathrm{e}^{-t_1x_1-t_2x_2}dx_1dx_2}{1+\mathrm{e}^{(b-x_1-x_2)/2}} = \frac{2\pi\mathrm{e}^{-t_2b}}{(t_1-t_2)\sin 2\pi t_2} -\sum_{k=1}^\infty\frac{(-1)^k\mathrm{e}^{-bk/2}}{(t_1-k/2)(t_2-k/2)}.
    \]
    The proof for the base case with $\mathrm{Re}(t_1)>\mathrm{Re}(t_2)$ follows similarly. Finally, if $\mathrm{Re}(t_1) = \mathrm{Re}(t_2)$ but $t_1\neq t_2$, the proof is completed by a continuity argument.

    If the conclusion holds for either $(n_1 - 1, n_2)$ or $(n_1, n_2-1)$, we prove it for $(n_1,n_2)$ similarly to the proof of \eqref{equ:3:2}. Finally, \eqref{equ:3:6} for $\mathrm{Re}(t_1)\in \frac{1}{2}\mathbb{Z}$ or $\mathrm{Re}(t_2)\in \frac{1}{2}\mathbb{Z}$ follows by continuity.
\end{proof}
\begin{proof}[Proof of Lemma \ref{lem:4:3}]
The convergence of the left-hand side is self-evident. We simplify it as follows:
\[
\begin{split}
    & \sum_{k=1}^\infty\frac{(-1)^k2t(2/k)^m}{t^2-k^2/4} = \sum_{k=1}^\infty(-1)^k(2/k)^m\left(\frac{1}{t - k/2} + \frac{1}{t + k/2}\right)\\
    & = \sum_{k\neq 0}\frac{(-1)^k(2/k)^m}{t-k/2} = \sum_{k\neq 0}\frac{(-1)^kt^{-m}}{t-k/2} + \sum_{k\neq 0}\frac{(-1)^k((2/k)^{m} - (1/t)^{m})}{t-k/2}.
\end{split}
\]
The first term simplifies as:
\[
\begin{split}
    & \sum_{k\neq 0}\frac{(-1)^kt^{-m}}{t-k/2} = t^{-m}\left(\sum_{k=-\infty}^\infty\frac{(-1)^k}{t-k/2} - \frac{1}{t}\right) \\
    & = t^{-m}\left(\frac{2\pi}{\sin 2\pi t}- \frac{1}{t}\right).
\end{split}
\]
The second term simplifies as:
\[
\begin{split}
    & \sum_{k\neq 0}\frac{(-1)^k((2/k)^{m} - (1/t)^{m})}{t-k/2} = \sum_{k\neq 0}\frac{(-1)^k(t^m-(k/2)^{m})}{t-k/2}(2/kt)^m \\
    & = \sum_{k\neq 0}(-1)^k\left(\sum_{i=0}^{m-1}(k/2)^{m-1-i}t^{i}\right)(2/kt)^m \\
    & = \sum_{i=0}^{m-1}\sum_{k\neq 0}(-1)^k(2/k)^{-1-i}t^{i-m} = \sum_{\substack{0<i<m\\ i\ odd}}\sum_{k\neq 0}(-1)^k(2/k)^{-1-i}t^{i-m}\\
    & = \sum_{\substack{0<i<m\\ i\ odd}}\frac{\sum_{k\neq 0}(-1)^k(2/k)^{-1-i}}{t^{m-i}},
\end{split}
\]
which is a Laurent series at $t=0$ with only the principal part. Since the left-hand side of \eqref{equ:4:3} is holomorphic, the second term should be exactly canceled with the principal part of the first term.
\end{proof}
\begin{proof}[Proof of Lemma \ref{lem:4:6}]
The equation is computed as follows:
\[
\begin{split}
    & \int_0^\infty ((x+\omega\mathrm{i})^m - (x-\omega\mathrm{i})^m)\mathrm{e}^{-tx}dx = \int_0^\infty \left(2\sum_{j=0}^{\lfloor(m-1)/2\rfloor} \binom{m}{2j+1}x^{m-2j-1}(\omega \mathrm{i})^{2j+1}\right)\mathrm{e}^{-tx}dx\\
    & = 2\mathrm{i}\sum_{j=0}^{\lfloor(m-1)/2\rfloor}(-1)^j\omega^{2j+1}\binom{m}{2j+1}\int_0^\infty x^{m-2j-1}\mathrm{e}^{-tx}dx\\
    & = 2\mathrm{i}\sum_{j=0}^{\lfloor(m-1)/2\rfloor}(-1)^j\omega^{2j+1}\binom{m}{2j+1}(m-2j-1)!t^{2j-m}= 2\mathrm{i}\sum_{2j+1\leq m} \frac{(-1)^j(\omega t)^{2j+1}}{(2j+1)!}\frac{m!}{t^{-m-1}} \\
    & = \mathsf{P}_t\left(2\mathrm{i}\sin \omega t \frac{m!}{t^{-m-1}} \right) = \mathsf{P}_t\left(2\mathrm{i}\sin \omega t\int_0^\infty x^m\mathrm{e}^{-tx}dx\right).
\end{split}
\]
\end{proof}